\documentclass[10pt]{amsart}

\input{diagrams}
\newtheorem{proposition}{Proposition}[section]
\newtheorem{lemma}[proposition]{Lemma}
\newtheorem{corollary}[proposition]{Corollary}
\newtheorem{theorem}[proposition]{Theorem}

\theoremstyle{definition}
\newtheorem{definition}[proposition]{Definition}

\theoremstyle{remark}
\newtheorem{remark}[proposition]{Remark}

\newcommand{\thlabel}[1]{\label{th:#1}}
\newcommand{\thref}[1]{Theorem~\ref{th:#1}}
\newcommand{\selabel}[1]{\label{se:#1}}
\newcommand{\seref}[1]{Section~\ref{se:#1}}
\newcommand{\lelabel}[1]{\label{le:#1}}

\newcommand{\prlabel}[1]{\label{pr:#1}}
\newcommand{\prref}[1]{Proposition~\ref{pr:#1}}
\newcommand{\colabel}[1]{\label{co:#1}}

\newcommand{\relabel}[1]{\label{re:#1}}

\newcommand{\delabel}[1]{\label{de:#1}}

\newcommand{\eqlabel}[1]{\label{eq:#1}}
\newcommand{\equref}[1]{(\ref{eq:#1})}

\newcommand{\Hom}{{\rm Hom}}

\newcommand{\ev}{{\rm ev}}
\newcommand{\coev}{{\rm coev}}

\def\lan{\langle}
\def\ran{\rangle}
\def\ot{\otimes}

\def\cal{\mathcal}

\newcommand{\Cc}{\mathcal{C}}
\newcommand{\Dd}{\mathcal{D}}

\newcommand{\Mm}{\mathcal{M}}

\newcommand{\Zz}{\mathcal{Z}}
\newcommand{\Ww}{\mathcal{W}}
\newcommand{\YD}{\mathcal{YD}}
\def\*C{{}^*\hspace*{-1pt}{\Cc}}
\def\text#1{{\rm {\rm #1}}}

\def\ol{\overline}

\def\yd{\mbox{$_H^H{\mathcal YD}$}}
\def\va{\varepsilon}
\def\v{\varphi}

\def\tr{\triangleright}
\def\rh{\rightharpoonup}
\def\lh{\leftharpoonup}

\def\ra{\rightarrow}
\def\a{\alpha}
\def\b{\beta}

\def\l{\lambda}
\def\r{\rho}
\def\cd{\cdot}

\def\ov{\overline}
\def\un{\underline}
\def\mf{\mathfrak}

\newcommand{\smi}{\mbox{$S^{-1}$}}

\def\rawo\lonra{\longrightarrow}

\def\ot{\otimes}

\newcommand{\tpla}{\mbox{$\tilde {p}^1$}}
\newcommand{\tplb}{\mbox{$\tilde {p}^2$}}

\newcommand{\tqla}{\mbox{$\tilde {q}^1$}}
\newcommand{\tqlb}{\mbox{$\tilde {q}^2$}}
\newcommand{\tQla}{\mbox{$\tilde {Q}^1$}}
\newcommand{\tQlb}{\mbox{$\tilde {Q}^2$}}

\begin{document}
\title[Yetter-Drinfeld categories for quasi-Hopf 
algebras]{Yetter-Drinfeld categories for quasi-Hopf algebras}
\author{D. Bulacu}
\address{Faculty of Mathematics, University of Bucharest,
RO-010014 Bucharest 1, Romania}
\email{dbulacu@al.math.unibuc.ro}

\author{S. Caenepeel}
\address{Faculty of Applied Sciences,
Vrije Universiteit Brussel, VUB, B-1050 Brussels, Belgium}
\email{scaenepe@vub.ac.be}
\urladdr{http://homepages.vub.ac.be/\~{}scaenepe/}
\author{F. Panaite}
\address{Institute of Mathematics of the Romanian
Academy, PO-Box 1-764, RO-70700 Bucha-rest, Romania}
\email{Florin.Panaite@imar.ro}

\thanks{Research supported by the bilateral project
"Hopf Algebras in Algebra, Topology, Geometry and Physics" of the 
Flemish and Romanian governments. This paper 
was finished while the first author 
was visiting the Vrije Universiteit Brussel, and he would like to
thank VUB for its warm hospitality. The third author was also partially 
supported by the programmes SCOPES and EURROMMAT}
\subjclass{16W30}

\keywords{quasi-Hopf algebra, Yetter-Drinfeld module, braided monoidal 
category}

\begin{abstract}
We show that all possible categories of Yetter-Drinfeld modules over a 
quasi-Hopf
algebra $H$ are isomorphic. We prove also that the category 
$\yd^{\rm fd}$ of finite dimensional left Yetter-Drinfeld 
modules is rigid and then we compute explicitly the 
canonical isomorphisms in $\yd^{\rm fd}$. Finally, we show that certain 
duals of $H_0$, the braided Hopf algebra introduced in 
\cite{bn,bpv}, are isomorphic as braided Hopf algebras if $H$ is a 
finite dimensional triangular quasi-Hopf algebra.
\end{abstract}
\maketitle

\section*{Introduction}
Let $H$ be a Hopf algebra with a bijective antipode. We can introduce
left, right, left-right and right-left Yetter-Drinfeld modules over $H$,
and it is well-known (see \cite{ag,rt}) that the corresponding 
categories $\yd $, ${\mathcal YD}^H_H$, ${}_H{\mathcal YD}^H$ and 
${}^H{\mathcal YD}_H$ are isomorphic. These categories are also 
isomorphic to the center of the monoidal category ${}_H\Mm$ of left 
$H$-modules, and, if $H$ is finite dimensional, to the category 
${}_{D(H)}\Mm$ of left modules over 
the Drinfeld double $D(H)$. It is also known that the category of finite 
dimensional Yetter-Drinfeld modules is rigid, that is, we have left and right 
duality in this category.

In \cite{dpr}, Dijkgraaf, Pasquier and Roche introduced the so-called
"twisted double" of a finite group, which is a Hopf algebra-type object
$D^{\omega}(G)$ associated to a pair $(G, \omega )$, where $G$ is a 
finite group and $\omega $ is a normalized $3$-cocycle on $G$; this object is
{\it not} a Hopf algebra, but a quasi-Hopf algebra in the sense of
Drinfeld \cite{d1}. The construction is similar to the quantum double, 
so it appears natural to try to define the quantum double of an arbitrary
finite dimensional quasi-Hopf algebra, generalizing the Drinfeld double
for Hopf algebras, and then to show that $D^{\omega}(G)$ is such a
quantum double. This has been done first by Majid in \cite{m1};
he first computed the center of the monoidal category
${}_H{\mathcal M}$ of left $H$-modules over the quasi-Hopf algebra $H$ 
(we will denote this center by $\yd$ and call its objects left Yetter-Drinfeld 
modules over $H$). Then he defined the quantum double $D(H)$ by an
implicit Tannaka-Krein reconstruction procedure, in such a way that
${}_{D(H)}{\mathcal M}\cong \yd$. An explicit construction of the
quantum double (as a so-called ``diagonal crossed product") has been 
given afterwards by Hausser and Nill in \cite{hn1,hn2}. They
identified the category of left modules over their quantum double with
the category ${}_{H^{\rm cop}}^{H^{\rm cop}}{\mathcal YD}$ (we will 
denote this category by ${}_H{\mathcal YD}^H$ and call 
its objects left-right Yetter-Drinfeld modules over $H$). 
Now the question arises whether $\yd \cong {}_H{\mathcal YD}^H$.

The first aim of this paper is to show that, indeed, 
the categories $\yd$ and ${}_H{\mathcal YD}^H$ 
(and also two other categories ${\mathcal YD}_H^H$ and 
${}^H{\mathcal YD}_H$ which we will introduce) are isomorphic, 
even in the situation where $H$ is not finite 
dimensional. In an earlier version of this paper, a computational proof of
this result was given, which was much more complicated than the corresponding
proof for Yetter-Drinfeld modules over coassociative Hopf algebras.
Viewing the categories of Yetter-Drinfeld modules as (left or right) 
centers of corresponding categories of modules, a more transparent 
approach is possible, and this is what we will
do in \seref{2}. This approach was suggested to us by the referee.
 
In \seref{3}, we show that the category $\yd^{\rm fd}$ of finite 
dimensional (left) Yetter-Drinfeld modules is rigid; the left and right 
duals are constructed explicitly. In an arbitrary rigid braided monoidal 
category $\Cc$, we have canonical isomorphisms 
$M\cong M^{**}$ and $(M\ot N)^*\cong M^*\ot N^*$. 
In \seref{4}, we compute these isomorphisms in the case $\Cc=\yd^{\rm fd}$. 
If we then specialize to finite dimensional left modules over a
quasitriangular quasi-Hopf algebra, we recover some results from 
\cite{bpv2}. 

Let $B$ be a braided Hopf algebra, that is a Hopf algebra in $\yd^{\rm 
fd}$. Then $B^*$ and ${}^*B$ are also braided Hopf algebras. 
In \seref{5}, we study the special case where
$B=H_0$, the braided Hopf algebra introduced in \cite{bn,bpv}. In 
particular, we prove that ${}^*H_0$ and $H^*_0$ are isomorphic braided Hopf algebras
if $H$ is a finite dimensional triangular 
quasi-Hopf algebra. They are also isomorphic to the coopposite of the 
braided Hopf algebra $\un{H}^*$ introduced in \cite{bc1}.

\section{Preliminary results}\selabel{1}
\subsection{Quasi-Hopf algebras}\selabel{1.1}
We work over a commutative field $k$. All algebras, linear spaces
etc. will be over $k$; unadorned $\ot $ means $\ot_k$. Following
Drinfeld \cite{d1}, a quasi-bialgebra is a fourtuple $(H, \Delta ,
\va , \Phi )$ where $H$ is an associative algebra with unit,
$\Phi$ is an invertible element in $H\ot H\ot H$, and $\Delta :\
H\ra H\ot H$ and $\va :\ H\ra k$ are algebra homomorphisms
satisfying the identities
\begin{eqnarray}
&&(id \ot \Delta )(\Delta (h))=
\Phi (\Delta \ot id)(\Delta (h))\Phi ^{-1},\label{q1}\\
&&(id \ot \va )(\Delta (h))=h,
\mbox{${\;\;\;}$}%
(\va \ot id)(\Delta (h))=h,\label{q2}
\end{eqnarray}
for all $h\in H$. $\Phi$ has to be a normalized $3$-cocycle, 
in the sense that 
\begin{eqnarray}
&&(1\ot \Phi)(id\ot \Delta \ot id) (\Phi)(\Phi \ot 1)= (id\ot id
\ot \Delta )(\Phi ) (\Delta \ot id \ot id)(\Phi ),\label{q3}\\
&&(id \ot \va \ot id )(\Phi )=1\ot 1.\label{q4}
\end{eqnarray}
The map $\Delta $ is called the coproduct or the 
comultiplication, $\va $ the counit and $\Phi $ the reassociator. 
As for Hopf algebras \cite{sw} we denote 
$\Delta (h)=h_1\ot h_2$ (summation understood), but since $\Delta$ is only 
quasi-coassociative we adopt the further convention  
$$
(\Delta \ot id)(\Delta (h))=h_{(1, 1)}\ot h_{(1, 2)}\ot h_2,
\mbox{${\;\;\;}$} (id\ot \Delta )(\Delta (h))=
h_1\ot h_{(2, 1)}\ot h_{(2,2)},
$$
for all $h\in H$. We will denote the tensor components of $\Phi$
by capital letters, and the ones of $\Phi^{-1}$ by small letters,
namely
\begin{eqnarray*}
&&\Phi=X^1\ot X^2\ot X^3=T^1\ot T^2\ot T^3=V^1\ot
V^2\ot V^3=\cdots\\
&&\Phi^{-1}=x^1\ot x^2\ot x^3=t^1\ot t^2\ot t^3=
v^1\ot v^2\ot v^3=\cdots
\end{eqnarray*}
$H$ is called a quasi-Hopf algebra if, moreover, there exists an
anti-automorphism $S$ of the algebra $H$ and elements $\a , \b \in
H$ such that, for all $h\in H$:
\begin{eqnarray}
&&S(h_1)\a h_2=\va (h)\a \mbox{${\;\;\;}$ and ${\;\;\;}$}
h_1\b S(h_2)=\va (h)\b ,\label{q5}\\
&&X^1\b S(X^2)\a X^3=1
\mbox{${\;\;\;}$ and${\;\;\;}$}
S(x^1)\a x^2\b S(x^3)=1.\label{q6}
\end{eqnarray}
Note that, in Drinfeld's original definition, the antipode of a quasi-Hopf algebra 
is required to be bijective. The axioms for a quasi-Hopf 
algebra imply that $\va (\a )\va (\b )=1$, so, 
by rescaling $\a $ and $\b $, we may assume without loss 
of generality that $\va (\a )=\va (\b )=1$ and $\va \circ S=\va $. 
The identities (\ref{q2}-\ref{q4}) also imply that 
\begin{equation}\label{q7}
(\va \ot id\ot id)(\Phi )= (id \ot id\ot \va )(\Phi )=1\ot 1.
\end{equation}
Together with a quasi-Hopf algebra
$H=(H, \Delta , \va , \Phi , S, \a , \b )$ we also have $H^{\rm op}$, 
$H^{\rm cop}$
and $H^{\rm op, cop}$ as quasi-Hopf algebras, where "op" means opposite
multiplication and "cop" means opposite comultiplication. The quasi-Hopf
structures are obtained by putting $\Phi_{\rm op}=\Phi^{-1}$,
$\Phi_{\rm cop}=(\Phi ^{-1})^{321}$, $\Phi _{\rm op, cop}=\Phi ^{321}$,
$S_{\rm op}=S_{\rm cop}=(S_{\rm op, cop})^{-1}=S^{-1}$, $\a _{\rm 
op}=\smi (\b )$,
$\b _{\rm op}=\smi (\a )$, $\a _{\rm cop}=\smi (\a )$, 
$\b _{\rm cop}=\smi (\b )$, $\a _{\rm op, cop}=\b $ and 
$\b _{\rm op, cop}=\a $.

Recall that the definition of a quasi-Hopf algebra is 
"twist covariant" in the following sense. An invertible element
$F\in H\ot H$ is called a {\sl gauge transformation} or {\sl
twist} if $(\va \ot id)(F)=(id\ot \va)(F)=1$. If $H$ is a
quasi-Hopf algebra and $F=F^1\ot F^2\in H\ot H$ is a gauge
transformation with inverse $F^{-1}=G^1\ot G^2$, then we can
define a new quasi-Hopf algebra $H_F$ by keeping the
multiplication, unit, counit and antipode of $H$ and replacing the
comultiplication, reassociator and the elements $\alpha$ and $\beta$
by
\begin{eqnarray}
&&\Delta _F(h)=F\Delta (h)F^{-1},\label{g1}\\
&&\Phi_F=(1\ot F)(id \ot \Delta )(F) \Phi (\Delta \ot id)
(F^{-1})(F^{-1}\ot 1),\label{g2}\\
&&\a_F=S(G^1)\a G^2,~~
\b_F=F^1\b S(F^2).\label{g3}
\end{eqnarray}
It is well-known that the antipode of a Hopf
algebra is an anti-coalgebra morphism. For a quasi-Hopf algebra,
we have the following statement: there exists a gauge
transformation $f\in H\ot H$ such that
\begin{equation} \label{ca}
f\Delta (S(h))f^{-1}= (S\ot S)(\Delta ^{\rm op}(h))
\mbox{,${\;\;\;}$for all $h\in H$,}
\end{equation}
where $\Delta ^{\rm op}(h)=h_2\ot h_1$. The element $f$ can be 
computed explicitly. First set
$$
A^1\ot A^2\ot A^3\ot A^4=(\Phi \ot 1) (\Delta \ot id\ot
id)(\Phi ^{-1}),
$$
$$
B^1\ot B^2\ot B^3\ot B^4=
(\Delta \ot id\ot id)(\Phi )(\Phi ^{-1}\ot 1).
$$
Then define $\gamma, \delta\in H\ot H$ by
\begin{equation} \label{gd}
\gamma =S(A^2)\a A^3\ot S(A^1)\a A^4~~{\rm and}~~ \delta
=B^1\b S(B^4)\ot B^2\b S(B^3).
\end{equation}
Then $f$ and $f^{-1}$ are given by the formulae:
\begin{eqnarray}
f&=&(S\ot S)(\Delta ^{\rm op}(x^1)) \gamma \Delta (x^2\b
S(x^3)),\label{f}\\
f^{-1}&=&\Delta (S(x^1)\a x^2) \delta (S\ot S)(\Delta
^{\rm op}(x^3)).\label{g}
\end{eqnarray}
Moreover, $f=f^1\ot f^2$ and $g=g^1 \ot g^2$ satisfy the  relations
\begin{equation} \label{gdf}
f\Delta (\a )=\gamma , ~~
\Delta (\b )f^{-1}=\delta 
\end{equation}
and (see \cite{bn3})
\begin{equation}\label{l3a}
g^1S(g^2\a)=\b ,~~
S(\b f^1)f^2=\a.
\end{equation}
Furthermore the corresponding twisted reassociator (see
(\ref{g2})) is given by
\begin{equation} \label{pf}
\Phi _f=(S\ot S\ot S)(X^3\ot X^2\ot X^1).
\end{equation}
In a Hopf algebra $H$, we obviously have the identity
$h_1\ot h_2S(h_3)=h\ot 1$, for all $h\in H$.
We will need the generalization of this formula to
the quasi-Hopf algebra setting.
Following \cite{hn1} and \cite{hn2}, we define
\begin{eqnarray}
&&p_R=p^1\ot p^2=x^1\ot x^2\b S(x^3),\label{qr}\\
&&q_R=q^1\ot q^2=X^1\ot S^{-1}(\a X^3)X^2,\label{qra}\\
&&p_L=\tilde{p}^1\ot \tilde{p}^2=X^2\smi (X^1\b )\ot X^3,\label{ql}\\
&&q_L=\tilde{q}^1\ot \tilde{q}^2=S(x^1)\a x^2\ot x^3.\label{qla}
\end{eqnarray}
For all $h\in H$, we then have:
\begin{eqnarray}
&&\Delta (h_1)p_R[1\ot S(h_2)]=p_R[h\ot 1],\label{qr1}\\
&&[1\ot S^{-1}(h_2)]q_R\Delta (h_1)=(h\ot 1)q_R,\label{qr1a}\\
&&[S(h_1)\ot 1]q_L\Delta (h_2)=(1\ot h)q_L,\label{ql1a}
\end{eqnarray}
and
\begin{eqnarray}
&&\Delta (q^1)p_R[1\ot S(q^2)]=1\ot 1,\label{pqr}\\
&&[S(\tilde{p}^1)\ot 1]q_L\Delta (\tilde{p}^2)=1\ot 1,\label{pql}\\
&&\Delta (\tilde{q}^2)p_L[\smi (\tilde{q}^1)\ot 1]=1\ot 1,\label{pqla}
\end{eqnarray}
\begin{eqnarray}
&&\hspace*{-2cm}
\Phi (\Delta \ot id_H)(p_R)
(p_R\ot id_H)\nonumber\\
&&=(id \ot \Delta )(\Delta (x ^1)p_R)
(1\ot g^1S(x ^3)\ot g^2S(x ^2)),\label{tpr2}
\end{eqnarray}
where $f=f^1\ot f^2$ is the twist defined in (\ref{f}) with its
inverse $f^{-1}=g^1\ot g^2$ defined in (\ref{g}).

\subsection{Quasitriangular quasi-Hopf algebras}\selabel{1.2}
A quasi-Hopf algebra $H$ is quasitriangular if
there exists an element $R\in H\ot H$ such that
\begin{eqnarray}
(\Delta \ot id)(R)&=&\Phi _{312}R_{13}\Phi ^{-1}_{132}R_{23}\Phi 
,\label{qt1}\\
(id \ot \Delta )(R)&=&\Phi ^{-1}_{231}R_{13}\Phi _{213}R_{12}\Phi 
^{-1},\label{qt2}\\
\Delta ^{\rm op}(h)R&=&R\Delta (h),~{\rm for~all~}h\in H,\label{qt3}\\
(\va \ot id)(R)&=&(id\ot \va)(R)=1.\label{qt4}
\end{eqnarray}
Here we use the following notation.
If $\sigma $ is a permutation of $\{1, 2, 3\}$, we set 
$\Phi _{\sigma (1)\sigma (2)\sigma (3)}=X^{\sigma ^{-1}(1)}\ot
X^{\sigma ^{-1}(2)}\ot X^{\sigma ^{-1}(3)}$, and
$R_{ij}$ means $R$ acting non-trivially
in the $i^{th}$ and $j^{th}$ positions of $H\ot H\ot H$.

In {\cite{bn3}} it is shown that $R$ is invertible. 
The inverse of $R$ is given by
\begin{equation}\label{invr1}
R^{-1}=X^1\b S(Y^2R^1x^1X^2)\a Y^3x^3X^3_2\ot Y^1R^2x^2X^3_1.
\end{equation}
Furthermore, the element
\begin{equation} \label{elmu}
u=S(R^2p^2)\a R^1p^1
\end{equation}
(with $p_R=p^1\ot p^2$ defined as in (\ref{qr})) 
satisfies $S^2(u)=u$, is invertible in $H$, and
\begin{equation} \label{inelmu}
u^{-1}=X^1R^2p^2S(S(X^2R^1p^1)\a X^3),
\end{equation}
\begin{equation} \label{sqina}
\va (u)=1~~{\rm and}~~
S^2(h)=uhu^{-1},
\end{equation}
for all $h\in H$. Consequently the antipode $S$ is bijective, so,
as in the Hopf algebra case, the assumptions about invertibility
of $R$ and bijectivity of $S$ can be dropped. Moreover, the 
$R$-matrix $R=R^1\ot R^2$ satisfies the identities (see
\cite{ac}, \cite{hn2}, \cite{bn3}):
\begin{eqnarray}
&&\hspace*{-2cm}f_{21}Rf^{-1}=(S\ot S)(R),\label{ext}\\
&&\hspace*{-2cm}S(R^2)\a R^1=S(\a )u,\label{sext}
\end{eqnarray}
where $f=f^1\ot f^2$ is the twist defined in (\ref{f}), and
$f_{21}=f^2\ot f^1$. In addition, a second formula for the inverse 
of $R$ is
\begin{equation}\label{invr2}
R^{-1}=\tilde{q}^2_1X^2R^1p^1\ot \tilde{q}^2_2X^3\smi 
(\tilde{q}^1X^1R^2p^2),
\end{equation}
where $p_R=p^1\ot p^2$ and $q_L=\tilde{q}^1\ot \tilde{q}^2$ 
are the elements defined by (\ref{qr}) and (\ref{qla}). Finally, 
recall that a quasi-Hopf algebra $(H, R)$
is called triangular if $R^{-1}=R_{21}$, where $R_{21}=R^2\ot R^1$.
\subsection{Monoidal categories}\selabel{1.3}
If $V$ is an object of a category $\Cc$, then the identity morphism $V\to V$
will also be denoted by $V$. The identity functor $\Cc\to \Cc$ will be denoted
by $\Cc$.

A monoidal category $\Cc=(\Cc,\ot, I,a,l,r)$ consists of a category $\Cc$, 
a functor $\ot:\ \Cc\times\Cc\to \Cc$, called the tensor product, an 
object $I\in \Cc$ called the unit object, and natural isomorphisms 
$a:\ \ot\circ (\ot\times \Cc)\to \ot\circ (\Cc\times \ot)$ (the associativity
constraint), $l:\ \ot\circ (I\times \Cc)\to \Cc$ (the left unit constraint) and 
$r:\ \ot\circ (\Cc\times I)\to \Cc$ (the right unit constraint). $a$ has to
satisfy the pentagon axiom, and $l$ and $r$ have to satisfy the triangle
axiom. We refer to \cite[XI.2]{k} for a detailed discussion. In the sequel,
we will identify $V\ot I\cong I\cong I\ot V$ using $l_V$ and $r_V$, for any
object $V\in \Cc$.

A monoidal functor between two monoidal categories $\Cc$ and 
$\Dd$ is a triple $(F,\varphi_0,\varphi_2)$, where $F:\ \Cc\to \Dd$ is a functor,
$\varphi_0:\ I\to F(I)$ is an isomorphism, and $\varphi_2(U,V):\
F(U)\ot F(V)\to F(U\ot V)$ is a family of natural isomorphisms in $\Dd$.
$\varphi_0$ and $\varphi_2$ have to satisfy certain properties, see for example
\cite[XI.4]{k}.

If $H$ is a quasi-bialgebra, then the categories ${}_H\Mm$ and $\Mm_H$ are
monoidal categories. The associativity constraint on ${}_H\Mm$ is the following:
for $U,V,W\in {}_H\Mm$, $a_{U,V,W}: (U\ot V)\ot W\to U\ot (V\ot W)$ is given by
\begin{equation}\eqlabel{1.3.1}
a_{U,V,W}((u\ot v)\ot w)=X^1\cd u\ot (X^2\cd v\ot X^3\cd w).
\end{equation}
On $\Mm_H$, the associativity constraint is given by the formula
\begin{equation}\eqlabel{1.3.2}
a_{U,V,W}((u\ot v)\ot w)=u\cd x^1\ot (v\cd x^2\ot w\cd x^3).
\end{equation}
Let $V\in \Cc$. $V^*\in \Cc$ is called a left dual of $V$, if the functor $-\ot V^*$
is the right dual of $-\ot V$. This is equivalent to the existence of morphisms
$\ev_V:\ V^*\ot V\to I$ and $\coev:\ I\to V\ot V^*$ such that
\begin{eqnarray}
&&(V\ot \ev_V)\circ a_{V,V^*,V}\circ (\coev_V\ot V)=V,\label{1.3.3}\\
&&(\ev_V\ot V^*)\circ a^{-1}_{V^*,V,V^*}\circ (V^*\ot \coev_V)=V^*,\label{1.3.4}
\end{eqnarray}
for all $V\in V$. ${}^*V\in \Cc$ is called a right dual of $V$ if
$-\ot {}^*V$ is the left dual of $-\ot V$. 
This is equivalent to the existence of morphisms 
$\ev'_V:\ V\ot {}^*V\to I$ and $\coev'_V:\ I\to {}^*V\ot V$ such that 
\begin{eqnarray}
&&(\ev'_V\ot V)\circ a^{-1}_{V,{}^*V,V}\circ (V^*\ot \coev'_V)=V,\label{1.3.5}\\
&&({}^*V\ot \ev'_V)\circ a_{{}^*V,V,{}^*V}\circ (\coev'_V\ot V^*)=V^*.\label{1.3.6} 
\end{eqnarray}
$\Cc$ is called a rigid monoidal category if every object of $\Cc$ has
a left and right dual. The category ${}_H\Mm^{\rm fd}$ of finite dimensional
modules over a quasi-Hopf algebra $H$ is rigid. For $V\in {}_H\Mm$,
$V^*=\Hom(V,k)$, with left $H$-action $\lan h\cdot \varphi,v\ran=
\lan\varphi,S(h)v\ran$. The evaluation and coevaluation are given by
\begin{eqnarray}
&&\hspace*{-2cm}
\ev_V(\varphi \ot v)=\varphi (\a \cd v),
\mbox{${\;\;\;}$}
\coev_V(1)=\b \cd v_i\ot v^i,\label{qrig}
\end{eqnarray}
where $\{v_i\}_i$ is a basis in
$V$ with dual basis $\{v^i\}_i$. The right dual ${}^*V$ of $V$ is
the same dual vector space now equipped with the left
$H$-module structure given  by $\lan h\cd \v ,v\ran=\lan\varphi,S^{-1}(h)v\ran$,
and with 
\begin{equation}
\ev'_V(v\ot \varphi )=\varphi (\smi (\a )\cd v),~~
\coev'_V(1)=v^i\ot \smi (\b )\cd v_i.\label{qrrig}
\end{equation}
The switch functor $\tau:\ \Cc \times\Cc\to \Cc \times \Cc$ 
is defined by $\tau(V,W)=(W,V)$.
A prebraiding on a monoidal category is a natural transformation 
$c:\ \ot\to\ot\circ\tau$, 
satisfying the hexagon axioms (see for example \cite[XIII.1]{k}). 
A prebraiding $c$ is called a braiding if it is a natural isomorphism. 
Let $(H,R)$ be a quasitriangular quasi-bialgebra. 
We then have the following  prebraiding $c$ on ${}_H\Mm$ 
(see \cite{k} or \cite{m2}):
\begin{equation}\label{br}
c_{U, V}(u\ot v)=R^2\cd v\ot R^1\cd u,
\end{equation}
which is a braiding if $H$ is a quasi-Hopf algebra. 

\subsection{Braided Hopf algebras}\selabel{1.4}
Let $\Cc$ be a braided monoidal category. We can define algebras, coalgebras,
bialgebras and Hopf algebras, extending the classical definitions from \cite{sw}
in the obvious way. 

Thus, a bialgebra in $\Cc $ is $(B, m, \eta , \Delta , \va )$ where
$B$ is an object in $\Cc $ and the morphism
$m:\  B\ot B\ra B$ gives a multiplication 
that is associative up to the isomorphism $a$. 
Similarly for the coassociativity of the comultiplication $\Delta : B\ra
B\ot B$. The identity in the algebra $B$ is expressed as usual
by $\eta : I\ra B$ such that 
$m\circ (\eta \ot id)=m\circ (id\ot \eta )=id$. The counit axiom is
$(\va \ot id)\circ \Delta =(id \ot \va )\circ \Delta =id$. 
In addition, $\Delta $ is required to be an algebra morphism where
$B\ot B$ has the multiplication $m_{B\ot B}$, 
defined as the composition 
\begin{equation}\label{bi}
\begin{matrix}
(B\ot B)\ot (B\ot B)\hspace*{-3mm}&\rTo^{a}& B\ot (B\ot (B\ot 
B))&\rTo^{B\ot a^{-1}}&
B\ot ((B\ot B)\ot B)\\
&\rTo^{B\ot c\ot B}&B\ot ((B\ot B)\ot B)&\rTo^{B\ot a}&B\ot (B\ot (B\ot 
B))\\
&\rTo^{a^{-1}}&(B\ot B)\ot (B\ot B)&\rTo^{m\ot m}& B\ot B.
\end{matrix}
\end{equation}
A Hopf algebra $B$ is a bialgebra with a
morphism $S: B\ra B$ in $\Cc$ (the antipode)
satisfying the usual axioms
$m\circ (S\ot id)\circ \Delta =\eta \circ \va =
m\circ (id \ot S)\circ \Delta $.

For a braided monoidal category ${\Cc}$, let ${\Cc}^{\rm in}$ be equal to $\Cc$
as a monoidal category, with the mirror-reversed braiding
$\tilde{c}_{M, N}=c^{-1}_{N, M}$. For a Hopf algebra $B\in \Cc$ with
bijective antipode, we  define $B^{\rm op}$, $B^{\rm cop}$ and $B^{\rm 
op, cop}$ by 
\begin{eqnarray}
&&m_{B^{\rm op}}=m_B\circ c^{-1}_{B, B},~~
\Delta _{B^{\rm op}}=\Delta _B,~~
S_{B^{\rm op}}=S^{-1}_B,\label{bop}\\
&&m_{B^{\rm cop}}=m_B,~~
\Delta _{B^{\rm cop}}=c^{-1}_{B, B}\circ \Delta _B,~~
S_{B^{\rm cop}}=S^{-1}_B,\label{bcop}\\
&&m_{B^{\rm op,cop}}=m_B\circ c_{B, B},~~
\Delta _{B^{\rm op,cop}}=c^{-1}_{B, B}\circ \Delta _B,~~
S_{B^{\rm op,cop}}=S_B,
\end{eqnarray}
and the other structure morphisms remain the same as for $B$. It
is well-known (see for example \cite{ag,ta}) 
that $B^{\rm op}$ and $B^{\rm cop}$ are Hopf algebras in ${\Cc}^{\rm in}$, 
and that $B^{\rm op,cop}$ is a Hopf algebra in $\Cc$.

\subsection{Monoidal categories and the center construction}\selabel{1.0}
Let $\Cc$ be a monoidal category. Following \cite{m1}, 
the weak left center $\Ww_l(\Cc)$ is the category
with the following objects and morphisms. An object is a couple $(V,s_{V, -})$,
with $V\in \Cc$ and $s_{V, -}:\ V\ot -\to -\ot V$ a natural transformation, 
satisfying the following condition, for all $X,Y\in \Cc$:
\begin{equation}\eqlabel{center1}
(X\ot s_{V, Y})\circ a_{X,V,Y}\circ (s_{V, X}\ot Y)=a_{X,Y,V}\circ s_{V, X\ot Y}
\circ a_{V,X,Y},
\end{equation}
and such that $s_{V, I}$ is the composition of the natural isomorphisms
$V\ot I\cong V\cong I\ot V$. A morphism between $(V,s_{V, -})$ and $(V',s_{V', -})$ 
consists of $\vartheta:\ V\to V'$ in $\Cc$ such that
$$
(X\ot \vartheta)\circ s_{V, X}=s_{V', X}\circ (\vartheta\ot X).
$$
$\Ww_l(\Cc)$ is a prebraided monoidal category. 
The tensor product is
$$
(V,s_{V, -})\ot (V',s_{V', -})=(V\ot V', s_{V\ot V', -})
$$
with
\begin{equation}\eqlabel{center2}
s_{V\ot V', X}=a_{X,V,V'}\circ (s_{V, X}\ot V')\circ a^{-1}_{V,X,V'}
\circ (V\ot s_{V', X})\circ a_{V,V',X},
\end{equation}
and the unit is $(I,I)$. The prebraiding $c$ on $\Ww_l(\Cc)$ is given by
\begin{equation}\eqlabel{center3}
c_{V,V'}=s_{V, V'}:\ (V, s_{V, -})\ot (V', s_{V', -})
\to (V',s_{V', -})\ot (V,s_{V, -}).
\end{equation}
The left center $\Zz_l(\Cc)$ is the full subcategory of $\Ww_l(\Cc)$
consisting of objects $(V, s_{V, -})$ with $s_{V, -}$ a natural isomorphism. 
$\Zz_l(\Cc)$ is a braided monoidal category. 
$\Zz_l(\Cc)^{\rm in}$ will be our notation for the
monoidal category $\Zz_l(\Cc)$, together with the inverse braiding
$\tilde{c}$ given by $\tilde{c}_{V,V'}=c^{-1}_{V',V}=s^{-1}_{V', V}$.

The weak right center $\Ww_r(\Cc)$ and the
right center $\Zz_r(\Cc)$ are introduced in a similar way. An object is a couple 
$(V,t_{-, V})$, where $V\in \Cc$ and
$t_{-, V}:\ -\ot V\to V\ot -$ is a family of natural transformations (resp. natural 
isomorphisms) such that $t_{-, I}$ is the natural isomorphism 
and 
\begin{equation}\eqlabel{center4}
a^{-1}_{V, X, Y}\circ t_{X\ot Y, V}\circ a^{-1}_{X, Y, V}=(t_{X, V}\ot Y)
\circ a^{-1}_{X, V, Y}\circ (X\ot t_{Y, V}),
\end{equation}
for all $X,Y\in \Cc$. 
A morphism between $(V,t_{-, V})$ and $(V',t_{-, V'})$ consists of 
$\vartheta:\ V\to V'$ in $\Cc$ such that
$$
(\vartheta\ot X)\circ t_{X, V}=t_{X, V'}\circ (X\ot \vartheta),
$$
for all $X\in \Cc$. $\Ww_r(\Cc)$ is prebraided monoidal and 
$\Zz_r(\Cc)$ is braided monoidal. The unit is $(I, I)$ and the 
tensor product is now
$$
(V,t_{-, V})\ot (V',t_{-, V'})=(V\ot V', t_{-, V\ot V'})
$$
with
\begin{equation}\eqlabel{center5}
t_{X, V\ot V'}=a^{-1}_{V,V',X}\circ (V\ot t_{X, V'})\circ a_{V,X,V'}
\circ (t_{X, V}\ot V')\circ a^{-1}_{X,V,V'}.
\end{equation}
The braiding $d$ is given by
\begin{equation}\eqlabel{center6}
d_{V,V'}=t_{V, V'}:\ (V, t_{-, V})\ot (V', t_{-, V'})\to 
(V', t_{-, V'})\ot (V, t_{-, V}).
\end{equation}
$\Zz_r(\Cc)^{\rm in}$ is the monoidal category $\Zz_r(\Cc)$ with the inverse
braiding $\tilde{d}$ given by $\tilde{d}_{V,V'}=d^{-1}_{V',V}=t^{-1}_{V', V}$.

For details in the case where $\Cc$ is a strict monoidal category, we refer
to \cite[Theorem XIII.4.2]{k}. The results remain valid in the case of
an arbitrary monoidal category, since every monoidal category is equivalent
to a strict one.

\begin{proposition}\prlabel{1.0.1}
Let $\Cc$ be a monoidal category. Then we have a braided isomorphism of
braided monoidal categories $F:\ \Zz_l(\Cc)\to \Zz_r(\Cc)^{\rm in}$, given by
$$
F(V, s_{V, -})=(V, s^{-1}_{V, -})~~{\rm and}~~F(\vartheta)=\vartheta .
$$
\end{proposition}

\begin{proof}
The proof is straighforward. Let us show that $F$ preserves the braiding.
Applying $F$ to the braiding map
$$
c_{V,V'}=s_{V, V'}:\ (V, s_{V, -})\ot (V',s_{V', -})\to 
(V', s_{V', -})\ot (V, s_{V, -}),
$$
we find
$$
F(s_{V, V'})=s_{V, V'}:\ (V, s^{-1}_{V, -})\ot (V', s^{-1}_{V', -})\to 
(V', s^{-1}_{V', -})\ot (V, s^{-1}_{V, -}).
$$
Write $t_{-, V}=s^{-1}_{V, -}$, $t_{-, V'}=s^{-1}_{V', -}$. With notation as above, 
we find
$$
\tilde{d}_{V,V'}=t^{-1}_{V', V}=s_{V, V'},
$$
as needed.
\end{proof}

Let $\Cc=(\Cc,\ot, I, a, l, r)$ be a monoidal category. Then we have
a second monoidal structure on $\Cc$, defined as follows:
$$
\ol{\Cc}=(\Cc,\ol{\ot}=\ot\circ\tau,I, \ol{a}, r, l)
$$
with $\tau:\ \Cc\times\Cc\to \Cc\times \Cc$, $\tau(V,W)=(W,V)$ and $\ol{a}$ 
defined by $\ol{a}_{V,W,X}=a^{-1}_{X,W,V}$. 

If $c$ is a (pre)braiding on $\Cc$, then $\ol{c}$, given by $\ol{c}_{V,W}=
c_{W,V}$ is a (pre)braiding on $\ol{\Cc}$. 

The following result is then completely obvious.

\begin{proposition}\prlabel{1.0.2}
Let $\Cc$ be a monoidal category. Then
$$
\ol{\Ww_l(\Cc)}\cong \Ww_r(\ol{\Cc})~~;~~
\ol{\Ww_r(\Cc)}\cong \Ww_l(\ol{\Cc}),
$$
as prebraided monoidal categories, and 
$$
\ol{\Zz_l(\Cc)}\cong \Zz_r(\ol{\Cc})~~;~~
\ol{\Zz_r(\Cc)}\cong \Zz_l(\ol{\Cc})
$$ 
as braided monoidal categories.
\end{proposition}

\begin{proposition}\prlabel{1.0.3}
Let $\Cc$ be a rigid monoidal category. Then the weak left 
(resp. right) center of $\Cc$
coincides with the  left 
(resp. right) center of $\Cc$, and is a rigid braided monoidal category.
\end{proposition}

\begin{proof}
For details in the case where $\Cc$ is strict, we refer to \cite[Lemma 7.2 and 7.3, Cor. 7.4]{s}. 
The general case then follows 
from the fact that every monoidal category is equivalent to a strict one.
For later use, we mention that for $(V,s)\in \Zz_l(\Cc)$, 
$(V, s_{V, -})^*=(V^*,s_{V^*, -})$, with
$s_{V^*, X}$ given by the following composition:
\begin{equation}\eqlabel{1.0.3.1}
\begin{matrix}
s_{V^*, X}:\ V^*\ot X&\rTo^{(V^*\ot X)\ot\coev_V}&(V^*\ot X)\ot(V\ot V^*)\\
&\rTo^{a_{V^*,X,V\ot V^*}}&V^*\ot (X\ot(V\ot V^*))\\
&\rTo^{V^*\ot a^{-1}_{X,V,V^*}}& V^*\ot ((X\ot V)\ot V^*)\\
&\rTo^{V^*\ot (s_{V, X}^{-1}\ot V^*)}& V^*\ot ((V\ot X)\ot V^*)\\
&\rTo^{V^*\ot a_{V,X,V^*}}& V^*\ot (V\ot (X\ot V^*))\\
&\rTo^{a^{-1}_{V^*,V,X\ot V^*}}& (V^*\ot V)\ot (X\ot V^*)\\
&\rTo^{\ev_V\ot X\ot V^*}& X\ot V^*.
\end{matrix}
\end{equation}
\end{proof}

\section{Yetter-Drinfeld modules over a quasi-Hopf algebra}\selabel{2}
\setcounter{equation}{0}
From \cite{m1}, we recall the notion of left Yetter-Drinfeld module over a
quasi-bialgebra. We also introduce right, left-right and right-left 
Yetter-Drinfeld modules. The aim of this Section is to study the relations 
between these four types of modules.

\begin{definition}\delabel{2.1}
Let $H$ be a quasi-bialgebra, with reassociator $\Phi$. 
A left Yetter-Drinfeld module 
is a left $H$-module $M$ together with a 
$k$-linear map (called the left $H$-coaction)
$$
\lambda_M:\ M\to H\ot M,~~\lambda_M(m)=m_{(-1)}\ot m_{(0)}
$$
such that the following conditions hold, for all $h\in H$ and $m\in M$:
\begin{eqnarray}
&&X^1m_{(-1)}\ot (X^2\cd m_{(0)})_{(-1)}X^3
\ot (X^2\cd m_{(0)})_{(0)}\nonumber\\
&&\hspace*{1cm}=X^1(Y^1\cd m)_{(-1)_1}Y^2\ot X^2(Y^1\cd 
m)_{(-1)_2}Y^3\ot X^3\cd (Y^1\cd m)_{(0)},\label{y1}\\
&&\va (m_{(-1)})m_{(0)}=m\label{y2},\\
&&h_1m_{(-1)}\ot h_2\cd m_{(0)}=(h_1\cd m)_{(-1)}h_2
\ot (h_1\cd m)_{(0)}.\label{y3}
\end{eqnarray}
The category of left Yetter-Drinfeld modules and $k$-linear
maps that preserve the $H$-action and $H$-coaction is denoted by
$\yd$.
\end{definition}

Let $H$ be a Hopf algebra. If $M$ is a left $H$-module and an 
$H$-comodule then (\ref{y3}) is equivalent to 
$$
(h\cd m)_{(-1)}\ot (h\cd m)_{(0)}=h_1m_{(-1)}S(h_3)\ot h_2\cd m_{(0)},
$$
for all $h\in H$ and $m\in M$. For quasi-Hopf algebras we have 
the following result. 

\begin{lemma}\lelabel{2.2}
Let $H$ be a quasi-Hopf algebra, $M\in{}_H\Mm$, and $\lambda:\
M\to H\ot M$ a $k$-linear map satisfying (\ref{y1}-\ref{y2}). Then (\ref{y3})
is equivalent to
\begin{equation}\label{y3p}
(h\cd m)_{(-1)}\ot (h\cd m)_{(0)}=
[q^1h_1]_1(p^1\cd m)_{(-1)}p^2S(q^2h_2)\ot [q^1h_1]_2\cd 
(p^1\cd m)_{(0)},
\end{equation}
for all $h\in H$, $m\in M$, where $p_R=p^1\ot p^2$ and
$q_R=q^1\ot q^2$ are the elements defined in (\ref{qr}-\ref{qra}).
\end{lemma}
\begin{proof}
Suppose that (\ref{y3}) holds. For any $h\in H$ and $m\in M$ we compute 
that
\begin{eqnarray*}
&&\hspace*{-2cm}(h\cd m)_{(- 1)}\ot (h\cd m)_{(0)}\\
&{{\rm (\ref{pqr})}\atop =}&(q^1_1p^1h\cd m)_{(- 1)}q^1_2p^2S(q^2)\ot
(q^1_1p^1h\cd m)_{(0)}\\
&{{\rm (\ref{qr1})}\atop =}&([q^1h_1]_1p^1\cd m)_{(- 1)}[q^1h_1]_2
p^2S(q^2h_2)\ot ([q^1h_1]_1p^1\cd m)_{(0)}\\
&{{\rm (\ref{y3})}\atop =}&
[q^1h_1]_1(p^1\cd m)_{(- 1)}p^2S(q^2h_2)\ot [q^1h_1]_2\cd 
(p^1\cd m)_{(0)}.
\end{eqnarray*}
Conversely, assume that (\ref{y3p}) holds; in particular
we have that 
\begin{equation}\label{syd}
\l _M(m)=q^1_1(p^1\cd m)_{(- 1)}p^2S(q^2)\ot q^1_2\cd 
(p^1\cd m)_{(0)},
\mbox{${\;\;}$$\forall $ $m\in M$.}
\end{equation}
We then compute, for all $h\in H$ and $m\in M$:
\begin{eqnarray*}
&&\hspace*{-2cm}(h_1\cd m)_{(- 1)}h_2\ot (h_1\cd m)_{(0)}\\
&{{\rm (\ref{y3p})}\atop =}&[q^1h_{(1, 1)}]_1(p^1\cd m)_{(- 1)}
p^2S(q^2h_{(1, 2)})h_2\ot [q^1h_{(1, 1)}]_2\cd (p^1\cd m)_{(0)}\\
&{{\rm (\ref{qr1a})}\atop =}&h_1q^1_1(p^1\cd m)_{(- 1)}p^2S(q^2)\ot
h_2q^1_2\cd (p^1\cd m)_{(0)}\\
&{{\rm (\ref{syd})}\atop =}&h_1m_{(- 1)}\ot h_2\cd m_{(0)}.
\end{eqnarray*}
\end{proof} 

From \cite{m1}, if $H$ is a quasi-bialgebra 
then the category of left Yetter-Drinfeld
modules is isomorphic to the weak left center of 
the monoidal category ${\cal C}=:{}_H\Mm$,  
$$
\Ww_l({}_H\Mm)\cong \yd . 
$$
The prebraided monoidal structure on 
$\Ww_l({}_H\Mm)$ induces a prebraided monoidal structure on 
$\yd$. This structure is such that the forgetful functor 
$\yd\to {}_H{\mathcal M}$ is monoidal. The coaction on the 
tensor product $M\ot N$ of two 
Yetter-Drinfeld modules $M$ and $N$ is given by 
\begin{eqnarray}
&&\hspace*{-2cm}
\lambda _{M\ot N}(m\ot n)=X^1(x^1Y^1\cd m)_
{(-1)} x^2(Y^2\cd n)_{(-1)}Y^3\nonumber\\
&&\hspace*{2cm}
\ot X^2\cd (x^1Y^1\cd m)_{(0)}\ot X^3x^3\cd 
(Y^2\cd n)_{(0)},\label{y4}
\end{eqnarray}
and the braiding is given by 
\begin{equation}\label{y5}
c_{M, N}(m\ot n)=m_{(-1)}\cd n\ot m_{(0)}.
\end{equation}
Moreover, if $H$ is a quasi-Hopf algebra then 
$\Ww_l({}_H\Mm)=\Zz_l({}_H{\cal M})\cong \yd $, and 
the inverse of the braiding is given by (see \cite{bn}):
\begin{equation}
c^{-1}_{M, N}(n\ot m)=\tilde{q}^2_1X^2\cd (p^1\cd
m)_{(0)} \ot S^{-1}(\tilde{q}^1X^1(p^1\cd
m)_{(-1)}p^2S(\tilde{q}^2_2X^3))\cd n\label{y6}
\end{equation}
where $p_R=p^1\ot p^2$ and $q_L=\tilde{q}^1\ot \tilde{q}^2$
are the elements defined by (\ref{qr}) and (\ref{qla}), respectively.

We also introduce left-right, right-left and right 
Yetter-Drinfeld modules. More explicitly: 

\begin{definition}
Let $H$ be a quasi-bialgebra, with reassociator $\Phi$.

1) A left-right Yetter-Drinfeld module is a left $H$-module $M$ together with a
$k$-linear map (called the right $H$-coaction)
$$
\lambda_M:\ M\to M\ot H,~~\rho_M(m)=m_{(0)}\ot m_{(1)}
$$
such that the following conditions hold, for all $h\in H$ and $m\in M$:
\begin{eqnarray}
&&\hspace*{-1cm}
(x^2\cdot m_{(0)})_{(0)}\otimes (x^2\cdot m_{(0)})_{(1)}x^1
\otimes x^3m_{(1)}\nonumber \\
&=&x^1\cdot (y^3\cdot m)_{(0)}\otimes x^2(y^3\cdot m)_{(1)_1}y^1
\otimes x^3(y^3\cdot m)_{(1)_2}y^2,\label{lry1}\\
&&\hspace*{-1cm}
\varepsilon (m_{(1)})m_{(0)}=m,\label{lry2}\\
&&\hspace*{-1cm}
h_1\cdot m_{(0)}\otimes h_2m_{(1)}=(h_2\cdot m)_{(0)}
\otimes (h_2\cdot m)_{(1)}h_1\label{lry3}.
\end{eqnarray}
The category of left-right Yetter-Drinfeld modules is denoted by 
${}_H{\mathcal YD}^H$.

2) A right-left Yetter-Drinfeld module is a right $H$-module $M$ together with a
$k$-linear map (called the left $H$-coaction)
$$
\lambda_M:\ M\to H\ot M,~~\lambda_M(m)=m_{(-1)}\ot m_{(0)}
$$
such that the following conditions hold, for all $h\in H$ and $m\in M$:
\begin{eqnarray}
&&
m_{(-1)}x^1\otimes x^3(m_{(0)}\cdot x^2)_{(-1)}\otimes
(m_{(0)}\cdot x^2)_{(0)}\nonumber\\
&&\hspace*{1cm}
=y^2(m\cdot y^1)_{(-1)_1}x^1\otimes y^3(m\cdot y^1)_{(-1)_2}x^2
\otimes (m\cdot y^1)_{(0)}\cdot x^3,\label{rly1}\\
&&
\varepsilon (m_{(-1)})m_{(0)}=m,\label{rly2}\\
&&
m_{(-1)}h_1\otimes m_{(0)}\cdot h_2=
h_2(m\cdot h_1)_{(-1)}\otimes (m\cdot h_1)_{(0)}\label{rly3}.
\end{eqnarray}
The category of right-left Yetter-Drinfeld modules is denoted by 
${}^H{\mathcal YD}_H$.

3) A right Yetter-Drinfeld module is a right $H$-module $M$ together with a
$k$-linear map (called the right $H$-coaction)
$$
\lambda_M:\ M\to M\ot H,~~\rho_M(m)=m_{(0)}\ot m_{(1)}
$$
such that the following conditions hold, for all $h\in H$ and $m\in M$:
\begin{eqnarray}
&&
(m_{(0)}\cd X^2)_{(0)}\otimes X^1(m_{(0)}\cdot X^2)_{(1)}\otimes
m_{(1)}X^3\nonumber\\
&&\hspace*{1cm}=(m\cdot Y^3)_{(0)}\cd X^1\otimes Y^1(m\cdot 
Y^3)_{(1)_1}X^2
\otimes Y^2(m\cdot Y^3)_{(1)_2}X^3,\label{ry1}\\
&&
\varepsilon (m_{(1)})m_{(0)}=m,\label{ry2}\\
&&
m_{(0)}\cd h_1\otimes m_{(1)}h_2=
(m\cdot h_2)_{(0)}\otimes h_1(m\cdot h_2)_{(1)}\label{ry3}.
\end{eqnarray}
The category of right Yetter-Drinfeld modules is denoted by 
${\mathcal YD}_H^H$.
\end{definition}

As in \cite{m1} the above definition of Yetter-Drinfeld
modules was given using the center construction. More precisely:

\begin{theorem}\thlabel{2.2}
Let $H$ be a quasi-bialgebra. Then we have the following isomorphisms of
categories:
$$
\Ww_r({}_H\Mm)\cong {}_H\YD^H~~;~~
\Ww_r(\Mm_H)\cong \YD^H_H~~;~~\Ww_l(\Mm_H)\cong {}^H\YD_H.
$$
If $H$ is a quasi-Hopf algebra, then these three weak centers are equal to
the centers.
\end{theorem}

\begin{proof} (sketch) Take $(M, t_{-, M})\in \Ww_r({}_H\Mm)$, and consider
$\r _M=(\eta _H\ot M)\circ t_{H, M}:\ M\to M\ot H$, that is,
$$
\r _M(m)=m_{(0)}\ot m_{(1)}=t_{H, M}(1\ot m).
$$
$\r _M$ determines $t_{-, M}$ completely: for $x\in X\in {}_H\Mm$, consider
the map
$$
f:\ H\to X,~~f(h)=h\cdot x
$$
in ${}_H\Mm$. The naturality of $t_{-, M}$ entails that
$$
(M\ot f)\circ t_{H, M}=t_{X, M}\circ (f\ot M),
$$
and therefore 
\begin{equation}\eqlabel{2.2.1}
t_{X, M}(x\ot m)=
t_{X, M}((f\ot M)(1\ot m))=m_{(0)}\ot f(m_{(1)})=m_{(0)}\ot m_{(1)}\cd x.
\end{equation}
In particular, $m=t_{k, M}(1\ot m)=m_{(0)}\ot m_{(1)}\cd 1_{k}=
\varepsilon(m_{(1)})m_{(0)}$, 
so (\ref{lry2}) holds. If we evaluate \equref{center4}, with $X=Y=H$, at
$1\ot 1\ot m$, we find (\ref{lry1}). Finally, 
$$
h\cdot t_{H, M}(1\ot m)=t_{H, M}(h_1\ot h_2\cd m)=
(h_2\cd m)_{(0)}\ot (h_2\cd m)_{(1)}h_1,
$$
and
$$
h\cdot t_{H, M}(1\ot m)=h\cdot (m_{(0)}\ot m_{(1)})=
h_1\cd m_{(0)}\ot h_2m_{(1)},
$$
proving (\ref{lry3}), and we have shown that
$(M,\r _M)$ is a left-right Yetter-Drinfeld module.

Conversely, if $(M,\r _M)$ is a left-right Yetter-Drinfeld module, 
$(M, t_{-, M})$, with $t_{-, M}$ given by \equref{2.2.1}, is an object of 
$\Ww_r({}_H\Mm)$. 
If $H$ has an antipode, then $t_{-, M}$ is invertible; 
the inverse is given by
\begin{equation}
t^{-1}_{N, M}(m\ot n)=q^1_1x^1S(q^2x^3(\tplb \cd m)_{(1)}\tpla )\cd n
\ot q^1_2x^2\cd (\tplb \cd m)_{(0)}\label{y6a}
\end{equation}
where $q_R=q^1\ot q^2$ and $p_L=\tilde{p}^1\ot \tilde{p}^2$
are the elements defined by (\ref{qra}) and (\ref{ql}), respectively.

The proof of the other two isomorphisms can be done in a similar way, we 
leave the details to the reader.
\end{proof}

The prebraided monoidal structure on $\Ww_r({}_H\Mm)$ induces a monoidal structure on
${}_H{\cal YD}^H$. This structure is such that the forgetful functor 
${}_H{\cal YD}^H\to {}_H{\mathcal M}$ is monoidal. 
Using \equref{1.3.1} and \equref{center5}, we find that the coaction on the tensor product 
$M\ot N$ of two left-right 
Yetter-Drinfeld modules $M$ and $N$ is given by: 
\begin{eqnarray}
&&\hspace*{-2cm}
{\r }_{M{\ot}N}(m{\ot} n)=x^1X^1\cd (y^2\cd m)_{(0)}
{\ot} x^2\cd (X^3y^3\cd n)_{(0)}\nonumber\\
&&\hspace*{2cm}
\ot x^3(X^3y^3\cd n)_{(1)}X^2(y^2\cd m)_{(1)}y^1,\label{slrms2}
\end{eqnarray}
for all $m\in M$, $n\in N$.
We have already observed that the functor forgetting the coaction is monoidal,
so
\begin{equation}
h\cd (m\ot n)=h_1\cd m\ot h_2\cd n.\label{slrms1}
\end{equation}
The braiding ${\mf{c}}$ can easily be deduced from \equref{center6} 
and \equref{2.2.1}:
${\mf{c}}_{M, N}:M{\ot}N\ra N{\ot}M$ is given by
\begin{equation}\label{slrbs}
{\mf{c}}_{M, N}(m{\ot}n)=n_{(0)}{\ot}n_{(1)}\cd m,
\end{equation}
for $m\in M$ and $n\in N$. In the case when $H$ is a 
quasi-Hopf algebra, the inverse of the braiding is given by
\begin{equation}\label{slribs}
{\mf{c}}_{M, N}^{-1}(n{\ot} m)=q^1_1x^1
S(q^2x^3(\tplb \cd n)_{(1)}\tpla )\cd m
\ot q^1_2x^2\cd (\tplb \cd n)_{(0)}.
\end{equation}
${{}_H{\mathcal YD}^H}^{\rm in}$ is the category with monoidal structure
(\ref{slrms2}-\ref{slrms1}) and the mirror reversed braiding 
$\tilde{{\mf c}}_{M,N}={\mf{c}}^{-1}_{N, M}$.

For completeness' sake, let us also describe the prebraided monoidal structure 
on $\YD_H^H$ and ${}^H\YD_H$. For
$M,N\in \YD_H^H$, the coaction on $M\ot N$ is given by the formula
\begin{eqnarray*}
&&\hspace*{-2cm}
\rho(m\ot n)=(m\cd X^2)_{(0)}\cd x^1Y^1 \ot (n\cd X^3x^3)_{(0)}\cd Y^2\\
&&\hspace*{2cm}
\ot X^1(m\cd X^2)_{(1)}x^2 (n\cd X^3x^3)_{(1)}Y^3,
\end{eqnarray*}
and the $H$-action is 
$$
(m\ot n)\cd h=m\cd h_1\ot n\cd h_2.
$$
The braiding is the following:
$$
d_{M,N}(m\ot n)=n_{(0)}\ot m\cd n_{(1)}.
$$
Now take $M,N\in {}^H\YD_H$. The coaction on $M\ot N$ is the following:
\begin{eqnarray*}
&&\hspace*{-2cm}
\lambda(m\ot n)=x^3(n\cd x^2)_{(-1)}X^2(m\cd x^1X^1)_{(-1)}y^1\\
&&\hspace*{2cm}
(m\cd x^1X^1)_{(0)}\cd y^2\ot (n\cd x^2)_{(0)}\cd X^3y^3,
\end{eqnarray*}
and the $H$-action is 
$$
(m\ot n)\cd h=m\cd h_1\ot n\cd h_2.
$$
The braiding is given by
$$
{\mf d}_{M,N}(m\ot n)=n\cd m_{(-1)}\ot m_{(0)}.
$$

\begin{remark}
Let $H$ be a quasi-bialgebra. If ${\cal C}={}_H{\cal M}$, one can 
easily check that $\ol{\cal C}\cong {}_{H^{cop}}{\cal M}$. Then, by 
\prref{1.0.2}, we have $\Ww_l(\ol{\cal C})\cong \ol{\Ww_r({\cal C}})$, 
that is ${}_{H^{cop}}^{H^{cop}}{\cal YD}\cong \ol{{}_H{\cal YD}^H}$ as 
prebraided monoidal categories; in an earlier version of this paper, we used 
this as a {\it definition} for left-right Yetter-Drinfeld modules in terms 
of left Yetter-Drinfeld modules.      
\end{remark}

Let $H$ be a quasi-bialgebra. We have that $H^{\rm op,cop}$ is also a
quasi-bialgebra.

\begin{proposition}\prlabel{2.3}
We have an isomorphism of monoidal categories
$$
F:\ \ol{{}_{H^{\rm op,cop}}{\Mm}}\to \Mm_H,
$$
given by $F(M)=M$ as a $k$-module, with right $H$-action $m\cdot h=h\cd m$.
\end{proposition}

\begin{proof}
It is well-known and obvious that $F$ is an isomorphism of categories.
So we only need to show that it preserves the monoidal structure. Let us
first describe the monoidal structure on ${}_{H^{\rm op,cop}}{\Mm}$.
The left $H^{\rm op,cop}$-action on $N\ot M$ is
$$
h\cd (n\ot m)=h_2\cd n\ot h_1\cd m.
$$
The associativity constraint $a_{P,N,M}:\ (P\ot N)\ot M\to P\ot (N\ot M)$
is
$$
a_{P,N,M}((p\ot n)\ot m)=X^3\cd p\ot(X^2\cd n\ot X^1\cd m).
$$
Now we describe the monoidal structure on $\ol{{}_{H^{\rm op,cop}}\Mm}$.
We have $M\ol{\ot}N=N\ot M$. For $m\in M$, $n\in N$, we write
$$
m\ol{\ot}n=n\ot m\in M\ol{\ot}N=N\ot M.
$$
Then
\begin{equation}\eqlabel{2.3.1}
h\cd (m\ol{\ot} n)=h_1\cd m\ol{\ot}h_2\cd n.
\end{equation}
The associativity constraint
$$
\ol{a}_{M,N,P}=a^{-1}_{P,N,M}:\
(M\ol{\ot}N)\ol{\ot}P=P\ot (N\ot M)\to M\ol{\ot}(N\ol{\ot}P)=(P\ot N)\ot M
$$
is given by
\begin{equation}\eqlabel{2.3.2}
\ol{a}_{M,N,P}((m\ol{\ot}n)\ol{\ot}p)=x^1\cd m\ol{\ot} 
(x^2\cd n\ol{\ot} x^3\cd p).
\end{equation}
It is then clear from (\ref{eq:2.3.1}-\ref{eq:2.3.2}) that $F$ preserves the
monoidal structure.
\end{proof}

An immediate consequence of \prref{1.0.2}, \thref{2.2} and \prref{2.3} is 
then the following.

\begin{proposition}\prlabel{2.4}
Let $H$ be a quasi-bialgebra. Then we have the following isomorphisms of
prebraided monoidal categories, induced by the functor $F$ from \prref{2.3}:
$$
\YD_H^H\cong \ol{{}_{H^{\rm op,cop}}^{H^{\rm op,cop}}\YD}~~{\rm and}~~
{}^H\YD_H\cong \ol{{}_{H^{\rm op,cop}}\YD^{H^{\rm op,cop}}}.
$$
\end{proposition}
\begin{proof}
We have the following isomorphisms of categories:
$$
\YD_H^H\cong \Ww_r(\Mm_H)\cong \Ww_r(\ol{{}_{H^{\rm op,cop}}\Mm})\cong 
\ol{\Ww_l({}_{H^{\rm op,cop}}\Mm)}\cong 
\ol{{}_{H^{\rm op,cop}}^{H^{\rm op,cop}}\YD}.
$$
\end{proof}

Combining \prref{1.0.1} and \thref{2.2}, we find the following result.

\begin{theorem}\thlabel{2.5}
Let $H$ be a quasi-Hopf algebra. Then we have an isomorphism of braided
monoidal categories
$$
F:\ {{}_H{\mathcal YD}^H}^{\rm in}\to \yd,
$$
defined as follows. For $M\in {}_H{\mathcal YD}^H$,
$F(M)=M$ as a left $H$-module; the left $H$-coaction is given by
\begin{equation}\label{s0}
\l _M(m)=m_{(-1)}\ot m_{(0)}=
q^1_1x^1S(q^2x^3(\tplb \cd m)_{(1)}\tpla )\ot q^1_2x^2\cd (\tplb 
\cd m)_{(0)},
\end{equation}
for all $m\in M$, where $q_R=q^1\ot q^2$ and $p_L=\tpla \ot 
\tplb $ are the 
elements defined by (\ref{qra}) and (\ref{ql}), and 
$\r _M(m)=m_{(0)}\ot m_{(1)}$ 
is the right coaction of $H$ on $M$. The functor $F$ sends a morphism 
to itself.
\end{theorem}
\begin{proof}
$F$ is nothing else then the composition of the isomorphisms
$$
{{}_H{\mathcal YD}^H}^{\rm in}\to \Zz_{r}({}_H\Mm)^{\rm in}
\to \Zz_{l}({}_H\Mm)\to \yd.
$$
For $M\in {{}_H{\cal YD}^H}^{\rm in}$, we compute that 
the corresponding left Yetter-Drinfeld
module structure $\lambda_M$ on $M$ is the following: 
\begin{eqnarray*}
\lambda_M(m)&=&t_{H, M}^{-1}(m\ot 1)=
\tilde{{\mf c}}_{M,H}(m\ot 1)={\mf{c}}^{-1}_{H, M}(m\ot 1)\\
&=&q^1_1x^1S(q^2x^3(\tplb \cd m)_{(1)}\tpla )\ot q^1_2x^2\cd (\tplb 
\cd m)_{(0)},
\end{eqnarray*}
as needed.
\end{proof}

In the same way, we have the following result.

\begin{theorem}\thlabel{2.6}
Let $H$ be a quasi-Hopf algebra. Then the categories $\YD{}^H_H$ and 
${}^H\YD_H^{\rm in}$ are isomorphic 
as braided monoidal categories.
\end{theorem}

We now recall some generalities about quasi-bialgebras. Let $H$ be a quasi-bialgebra,
and ${\mf F}={\mf F}^1\ot {\mf F}^2\in H\ot H$ a twist with inverse
${\mf F}^{-1}={\mf G}^1\ot {\mf G}^2$. Then we have an isomorphism of
monoidal categories
$$
\Pi:\ {}_H\Mm\to {}_{H_{\mf F}}\Mm.
$$
$\Pi(M)=M$, with the same left $H$-action. If $H$ is a quasi-Hopf algebra, then we can
consider the Drinfeld twist $f$ defined in (\ref{f}). The antipode
$S:\ H^{\rm op,cop}\to H_f$ is a quasi-Hopf algebra isomorphism, and therefore
the monoidal categories ${}_{H^{\rm op,cop}}\Mm$ and ${}_{H_f}\Mm$ are isomorphic. 
We have seen in \prref{2.3} that ${}_{H^{\rm op,cop}}\Mm$ is isomorphic to
$\ol{\Mm _H}$ as a monoidal category.
We conclude that the monoidal categories ${}_H\Mm$ and $\ol{\Mm _H}$ are isomorphic.
Using \prref{1.0.2} and \thref{2.2}, we find braided monoidal isomorphisms
$$
{}_H^H\YD\cong \Zz_l({}_H\Mm)\cong \Zz_l(\ol{\Mm _H})\cong
\ol{\Zz_r({\Mm}_H)}\cong \ol{\YD _H^H}
$$
and
$$
{}_H\YD^H\cong \Zz_r({}_H\Mm)\cong \Zz_r(\ol{\Mm _H})\cong
\ol{\Zz_l(\Mm_H)}\cong \ol{{}^H\YD _H}.
$$

We summarize our results as follows:

\begin{theorem}\thlabel{2.7}
Let $H$ be a quasi-Hopf algebra. Then we have the following isomorphisms
of braided monoidal categories:
$$
{}_H^H\YD\cong {{}_H\YD^H}^{\rm in}\cong \ol{\YD _H^H}
\cong \ol{{}^H\YD_H}^{\rm in}.$$
\end{theorem}

The isomorphisms ${}_H^H\YD\cong \ol{\YD _H^H}$ and 
${}_H\YD^H\cong \ol{{}^H\YD_H}$ 
can be described explicitely. Let us compute the functor 
${}_H^H\YD\to \ol{\YD _H^H}$.

We have a monoidal isomorphism $\Pi:\ {}_H\Mm\to {}_{H_f}\Mm$; $\Pi(M)=M$ with
the same left $H$-action. Denote the tensor product on ${}_{H_f}\Mm$ by $\ot^f$.
For $M,N\in {}_H\Mm$, the isomorphism $\psi:\ \Pi(M\ot N)\to \Pi(M)\ot ^f\Pi(N)$
is given by
$$
\psi(m\ot n)=f^1\cd m\ot ^ff^2\cd n.
$$
This isomorphism induces an isomorphism between the two
left centers, hence between the categories ${}_H^H\YD$ and ${}_{H_f}^{H_f}{\YD}$.
Take $(M, s_{M, -})\in \Zz_l({}_H\Mm)$ and $(M, s^f_{M, -})$ the corresponding 
object in $\Zz_l({}_{H_f}\Mm)$, and let 
$\lambda$ and $\lambda^f$ be the associated 
coactions. We then have the following commutative diagram
$$
\begin{diagram}
M\ot H&\rTo^{s_{M, H}}& H\ot M\\
\dTo^{\psi}&&\dTo_{\psi}\\
M\ot^f H&\rTo^{s^f_{M, H}}& H\ot^f M
\end{diagram}
$$
and we compute that
$$
\lambda^f(m)=s_{M, H}^f(m\ot ^f1)=\psi((s_{M, H}(g^1\cd m\ot g^2))=
f^1(g^1\cd m)_{(-1)}g^2\ot ^ff^2(g^1\cd m)_{(0)}.
$$
The quasi-Hopf algebra isomorphism $S^{-1}:\ H_f\to H^{\rm op,cop}$
induces an isomorphism of monoidal categories
$$
{}_{H_f}\Mm\to {}_{H^{\rm op,cop}}\Mm\cong \ol{\Mm _H}.
$$
The image of $M$ is $M$ as a $k$-vector space, with right $H$-action given by 
$m\cd h=S(h)\cd m$. Take $(M, s^f_{M, -})\in \Zz_l({}_{H_f}\Mm)$ and the
corresponding object $(M, t_{-, M})\in \ol{\Zz_r(\Mm_H)}$. 
Using similar arguments as above we can compute the right $H$-coaction 
$\rho$ on $M$:
$$
\rho _M(m)=f^2\cd (g^1\cd m)_{(0)}\ot 
S^{-1}(f^1(g^1\cd m)_{(-1)})g^2).
$$
We conclude that the braided isomorphism $K:\ \yd\to \ol{\YD_H^H}$ 
is defined as follows: $K(M)=M$, with
\begin{eqnarray*}
&&m\cd h=S(h)\cd m,\\
&&
\r _M(m)=f^2\cd (g^1\cd m)_{(0)}\ot \smi (f^1(g^1\cd m)_{(-1)}g^2).
\end{eqnarray*}
The inverse functor $K^{-1}$ can be computed in a similar way: $K^{-1}(M)=M$
with
\begin{eqnarray*}
&&h\cd m=m\cd \smi (h),\\
&&\l _M(m)=g^1S((m\cd \smi (f^1))_{(1)})f^2\ot (m\cd \smi 
(f^1))_{(0)}\cd \smi (g^2).
\end{eqnarray*}
For completeness' sake, we also give the formulas for the 
braided monoidal isomorphism
$G:\ \ol{{}^H\YD_H}\to {}_H\YD^H$. $G(M)=M$ with
\begin{eqnarray*}
&&h\cd m=m\cd S^{-1}(h),\\
&&\r _M(m)=g^1\cd (f^2\cd m)_{(0)}\ot g^2S((f^2\cd m)_{(-1)})f^1.
\end{eqnarray*}
Conversely, $G^{-1}(M)=M$ with
\begin{eqnarray*}
&&m\cd h=S(h)\cd m,\\
&&
\l _M(m)=\smi (f^2(g^2\cd m)_{(1)}g^1)\ot f^1\cd (g^2\cd m)_{(0)}.
\end{eqnarray*}
\section{The rigid braided category $\yd ^{\rm fd }$}\selabel{3}
\setcounter{equation}{0}
It is well-known that the category of finite dimensional Yetter-Drinfeld
modules over a coassociative Hopf algebra with invertible antipode is
rigid. By \prref{1.0.3} and since ${}_H{\cal M}^{fd}$ is rigid the same 
result holds for the 
category of finite dimensional Yetter-Drinfeld
modules over a quasi-Hopf algebra. We will give the explicit formulas in
this Section. We first need a lemma.

\begin{lemma}\lelabel{3.1}
Let $H$ be a quasi-Hopf algebra and $p_R=p^1\ot p^2$, 
$q_L=\tqla \ot \tqlb =\tQla \ot \tQlb $ and $f=f^1\ot f^2$ the elements 
defined by (\ref{qr}), (\ref{qla}) and (\ref{f}), respectively. Then the 
following relations hold:
\begin{eqnarray}
&&\tilde{q}^1X^1\ot \tilde{q}^2_1X^2\ot \tilde{q}^2_2X^3=
S(x^1)\tqla x^2_1\ot \tqlb x^2_2\ot x^3,\label{fo1}\\
&&S(p^1)\tqla p^2_1S(\tQlb )_1\ot \tQla \tqlb p^2_2S(\tQlb )_2=f.\label{fo2}
\end{eqnarray}  
\end{lemma}
\begin{proof}
The equality (\ref{fo1}) follows easily from (\ref{q3}), (\ref{q5}) and the 
definition of $q_L$. In order to prove (\ref{fo2}), 
we denote by $\delta =\delta ^1\ot \delta ^2$ the element defined in (\ref{gd}), 
and then compute that 
\begin{eqnarray*}
&&\hspace*{-2cm}
S(p^1)\tqla p^2_1S(\tQlb )_1\ot \tQla \tqlb p^2_2S(\tQlb )_2\\
&{{\rm (\ref{qr})}\atop =}&
S(x^1)\tqla x^2_1\b _1S(\tQlb x^3)_1\ot \tQla \tqlb x^2_2\b _2S(\tQlb x^3)_2\\
&{{\rm (\ref{ca}, \ref{gdf})}\atop =}&
S(x^1)\tqla x^2_1\delta ^1S(\tilde{Q}^2_2x^3_2)f^1\ot 
\tQla \tqlb x^2_2\delta ^2S(\tilde{Q}^2_1x^3_1)f^2\\
&{{\rm (\ref{gd}, \ref{q3}, \ref{q5})}\atop =}&
S(z^1x^1)\a z^2x^2_1y^1\b S(\tilde{Q}^2_2x^3_2y^3_2X^3)f^1\ot 
\tQla z^3x^2_2y^2X^1\b S(\tilde{Q}^2_1x^3_1y^3_1X^2)f^2\\
&{{\rm (\ref{q3}, \ref{q5})}\atop =}&
S(x^1)\a x^2\b S(\tilde{Q}^2_2X^3x^3)f^1\ot 
\tQla X^1\b S(\tilde{Q}^2_1X^2)f^2\\
&{{\rm (\ref{q6}, \ref{ql})}\atop =}&
S(\tilde{Q}^2_2\tplb )f^1\ot \tQla S(\tilde{Q}^2_1\tpla )f^1\\
&{{\rm (\ref{pqla})}\atop =}&f^1\ot f^2=f,
\end{eqnarray*} 
as needed, and this finishes the proof. 
\end{proof}

\begin{theorem}\thlabel{3.2}
Let $H$ be a quasi-Hopf algebra. Then $\yd ^{\rm fd}$ is a
braided monoidal rigid category. For a finite dimensional left Yetter-Drinfeld module
$M$ with basis $({}_im)_{i=\ov {1, n}}$ and corresponding dual basis 
$({}^im)_{i=\ov {1, n}}$, the left and right duals $M^*$ and ${}^*M$ are equal to
$\Hom(M,k)$ as a vector space, with the following $H$-action and $H$-coaction:
\begin{itemize}
\item[-]For $M^*$:
\begin{eqnarray}
&&(h\cd m^*)(m)=m^*(S(h)\cd m),\label{rdy1}\\
&&\l _{M^*}(m^*)=m^*_{(- 1)}\ot m^*_{(0)}\nonumber\\
&&\hspace*{1cm}
=\langle m^*, f^2\cd (g^1\cd {}_im)_{(0)}\rangle 
\smi (f^1(g^1\cd {}_im)_{(- 1)}g^2)\ot {}^im\label{rdy2}.
\end{eqnarray}
\item[-]For ${}^*M$:
\begin{eqnarray}
&&(h\cd {}^*m)(m)={}^*m(\smi (h)\cd m),\label{ldy1}\\
&&\l _{{}^*M}({}^*m)={}^*m_{(- 1)}\ot {}^*m_{(0)}\nonumber\\
&&\hspace*{1cm}
=\langle {}^*m,
\smi (f^1)\cd (\smi (g^2)\cd {}_im)_{(0)}\rangle
g^1S((\smi (g^2)\cd {}_im)_{(- 1)})f^2\ot {}^im,\label{ldy2}
\end{eqnarray}
\end{itemize}
for all $h\in H$, $m^*\in M^*$, ${}^*m\in {}^*M$ and $m\in M$. 
Here $f=f^1\ot f^2$ is the twist defined by
(\ref{f}), with inverse $f^{-1}=g^1\ot g^2$.
\end{theorem}

\begin{proof}
The left $H$-action on $M^*$ viewed as an element of $\yd$ is the same as the left
$H$-action on $M^*$ viewed as an element of ${}_H\Mm$.
We compute the left $H$-coaction, using \equref{1.0.3.1}. By (\ref{y6}) 
in $\Zz _l({}_H{\cal M})\cong \yd$ we have 
$$
s^{-1}_{V, X}(v\ot x)
=\tilde{q}^2_1X^2\cd (p^1\cd
x)_{(0)} \ot S^{-1}(\tilde{q}^1X^1(p^1\cd
x)_{(-1)}p^2S(\tilde{q}^2_2X^3))\cd v
$$
for all $V, X\in {}_H{\cal M}$, $x\in X$ and $v\in V$. Now, if we 
denote by $P^1\ot P^2$ another copy of $p_R$ 
and by $\tQla \ot \tQlb$ another copy of $q_L$ we then find
\begin{eqnarray*}
&&\hspace*{-15mm} \lambda(m^*)=s_{M^*, H}(m^*\ot 1)
=\lan x^1Z^1\cd m^*,\alpha x^2Y^1\tilde{q}^2_1X^2\cd 
\bigl(p^1y^2Z^3_1\beta \cd {}_im\bigr)_{(0)}\ran\\
&&x_1^3Y^2\tilde{q}^2_2X^3S^{-1}\Bigl(\tilde{q}^1X^1\bigl
(p^1y^2Z_1^3\beta \cd {}_im\bigr)_{(-1)}
p^2\Bigr)y^1Z^2\ot x_2^3Y^3y^3Z_2^3\cd {}^im\\
&{{\rm (\ref{q5}, \ref{qr})}\atop =}&
\lan m^*, \tQla Y^1\tilde{q}^2_1X^2\cd 
\bigl(p^1P^2S(\tilde{Q}^2_2Y^3)\cd {}_im\bigr)_{(0)}\ran\\
&&\tilde{Q}^2_1Y^2\tilde{q}^2_2X^3\smi \Bigl(\tilde{q}^1X^1
\bigl(p^1P^2S(\tilde{Q}^2_2Y^3)\cd {}_im\bigr)_{(-1)}p^2\Bigr)P^1\ot {}^im\\
&{{\rm (\ref{fo1}, \ref{y3})}\atop =}&
\lan m^*, \tQla Y^1\tqlb \cd \bigl(y^2_1p^1P^2S(\tilde{Q}^2_2Y^3)
\cd {}_im\bigr)_{(0)}\ran \\
&&\tilde{Q}^2_1Y^2y^3\smi \Bigl(\tqla \bigl(y^2_1p^1P^2S(\tilde{Q}^2_2Y^3)
\cd {}_im\bigr)_{(-1)}y^2_2p^2\Bigr)y^1P^1\ot {}^im\\
&{{\rm (\ref{tpr2}, \ref{y3})}\atop =}&
\lan m^*, \tQla Y^1\tqlb x^1_{(2, 2)}p^2_2\cd \bigl(g^1
S(\tilde{Q}^2_2Y^3x^3)\cd {}_im\bigr)_{(0)}\ran \\
&&\tilde{Q}^2_1Y^2\smi \Bigl(\tqla x^1_{(2, 1)}p^2_1\cd 
\bigl(g^1S(\tilde{Q}^2_2Y^3x^3)\cd {}_im\bigr)_{(-1)}g^2S(x^2)\Bigr)
x^1_1p^1\ot {}^im\\
&{{\rm (\ref{ql1a}, \ref{ca})}\atop =}&
\lan m^*, \tQla \tqlb p^2_2\cd \bigl(S(\tQlb )_1g^1\cd {}_im\bigr)_{(0)}\ran \\
&&\smi \Bigl(\tqla p^2_1\cd \bigl(S(\tQlb )_1g^1\cd {}_im\bigr)_{(-1)}
S(\tQlb )_2g^2\Bigr)p^1\ot {}^im\\
&{{\rm (\ref{y3}, \ref{fo2})}\atop =}&
\lan m^*, f^2\cd \bigl(g^1\cd {}_im\bigr)_{(0)}\ran 
\smi \Bigl(f^1\cd \bigl(g^1\cd {}_im\bigr)_{(-1)}g^2\Bigr)\ot {}^im,
\end{eqnarray*}
as claimed. The structure on ${}^*M$ can be computed in a similar way, 
we leave the details to the reader.
\end{proof}

\section{The canonical isomorphisms in $\yd^{\rm fd}$}\selabel{4}
\setcounter{equation}{0}
If ${\mathcal C}$ is a braided rigid category then, for any two
objects $M, N\in {\mathcal C}$, there exist two canonical
isomorphisms in ${\mathcal C}$: $M\cong M^{**}$ and $(M\ot N)^*\cong 
M^*\ot N^*$.
In this section our goal is to compute explicitly the above 
isomorphisms in the
particular case ${\mathcal C}=\yd ^{f. d.}$. Then we will specialize 
them for the
category of finite dimensional modules over a quasitriangular 
quasi-Hopf algebra.\\
Let ${\mathcal C}$ be a rigid monoidal category and $M, N$ objects of 
${\mathcal C}$.
If $\nu : M\ra N$ is a morphism in ${\mathcal C}$,
following \cite{k} we can define the transposes of $\nu $ as being the 
compositions; the left and right unit constraints are treated as the identity.
\begin{equation}
\begin{matrix}
\nu^*: \ N^*&
\rTo^{{N^*}\ot \coev_M}&N^*\ot (M\ot M^*)
&\rTo^{{N^*}\ot (\nu \ot id_{M^*})}&N^*\ot (N\ot M^*)\\
&
\rTo^{a^{-1}_{N^*, N, M^*}}&(N^*\ot N)\ot M^*
&\rTo^{\ev_N\ot {M^*}}&M^*,
\end{matrix}
\label{rt}
\end{equation}
and
\begin{equation}
\begin{matrix}
{}^*\nu : {}^*N&
\rTo^{\coev'_M\ot {{}^*N}}&({}^*M\ot M)\ot {}^*N
&\rTo^{({{}^*M}\ot \nu )\ot {{}^*N}}&
({}^*M\ot N)\ot {}^*N\\
&
\rTo^{a_{{}^*M, N, {}^*N}}&{}^*M\ot (N\ot {}^*N)
&\rTo^{{{}^*M}\ot \ev'_N}&{}^*M.
\end{matrix}\label{lt}
\end{equation}
Since the functors $-\ot V$ and $-\ot {}^*(V^*)$ are left duals of $-\ot V^*$,
they are naturally isomorphic, so we have
an isomorphism $\theta_M:\ V\cong {}^*(V^*)$, and this isomorphism is
natural in $M$. In the same way, we have a natural isomorphism 
$\theta'_M:\ V\cong ({}^*V)^*$.
Let us describe $\theta_M$ and its invers explicitly, we refer
to \cite{m2} for details.
\begin{equation}
\begin{matrix}
\theta _M:\ M& \rTo^{M\ot \coev_{{}^*M}}&M\ot ({}^*M\ot ({}^*M)^*)&&\\
&\rTo^{a^{-1}_{M, {}^*M, ({}^*M)^*}}
&
(M\ot {}^*{M})\ot ({}^*M)^*
&
\rTo^{\ev'_M\ot {({}^*M)^*}}&({}^*M)^*;
\end{matrix}\label{theta}
\end{equation}
\begin{equation}
\begin{matrix}
{\theta}_M^{-1}:\ ({}^*M)^*&
\rTo^{{({}^*M)^*}\ot \coev'_M}&({}^*M)^*\ot ({}^*M\ot M)&&\\
&\rTo^{a^{-1}_{({}^*M)^*, {}^*M, M}}&
(({}^*M)^*\ot {}^*M)\ot M
&\rTo^{\ev_{{}^*M}\ot M}&M.
\end{matrix}\label{thetam}
\end{equation}
We also have a natural isomorphism $\Theta _M:\ M^*\to {}^*M$. this can
be described as follows, see \cite{ag} for details.
\begin{equation}
\begin{matrix}
\Theta _M:\ M^*&
\rTo^{{M^*}\ot \coev'_M}&M^*\ot ({}^*M\ot M)
&\rTo^{a^{-1}_{M^*, {}^*M, M}}&(M^*\ot {}^*M)\ot M\\
&
\rTo^{c_{M^*, {}^*M}\ot M}&({}^*M\ot M^*)\ot M
&\rTo^{a_{{}^*M, M^*, M}}&{}^*M\ot (M^*\ot M)\\
&\rTo^{{{}^*M}\ot \ev_M}&{}^*M;
\end{matrix}\label{rl}
\end{equation}
\begin{equation}
\begin{matrix}
\Theta ^{-1}_M:\ {}^*M
&\rTo^{\coev_M\ot id_{{}^*M}}&(M\ot M^*)\ot {}^*M
&\rTo^{a_{M, M^*, {}^*M}}&M\ot (M^*\ot {}^*M)\\
&\rTo^{M\ot c^{-1}_{M^*, {}^*M}}&M\ot ({}^*M\ot M^*)
&\rTo^{a^{-1}_{M, {}^*M, M^*}}&(M\ot {}^*M)\ot M^*\\
&\rTo^{\ev'_M\ot {M^*}}&M^*.
\end{matrix}\label{lr}
\end{equation}
Thus the functors $(-)^*$ and ${}^*(-)$ are
naturally isomorphic, and we conclude that
$$
M^{**}=(M^*)^*\cong {}^*(M^*)\cong M\cong ({}^*M)^*\cong 
{}^*({}^*M)={}^{**}M.
$$
We will apply these results in the particular case where
${\mathcal C}=\yd^{\rm fd}$.

1) The maps $\nu ^*$ and ${}^*\nu $ coincide
with the usual transposed map of $\nu $. Indeed, by (\ref{rt}) and 
(\ref{q6}), we have that
\begin{eqnarray*}
\nu ^{*}(n^*)&=&\langle x^1\cd n^*, \a x^2\b \cd 
\nu ({}_in)\rangle x^3\cd {}^in\\
&=&\langle n^*, S(x^1)\a x^2\b S(x^3)\cd \nu 
({}_in)\rangle {}^in=n^*\circ \nu ,
\end{eqnarray*}
where $({}_in)_{i=\ov {1, t}}$ is a basis of $N$ and $({}^in)_{i=\ov 
{1, t}}$ its dual basis. A similar computation shows that 
${}^*\nu ({}^*n)={}^*n\circ \nu $ for any ${}^*n\in {}^*N$.

2) It is not hard to see that
the map $\theta _M$ defined by (\ref{theta}) is given by
\begin{equation}\label{thetay}
\theta _M(m)=\langle {}^im, m\rangle {}^{*i}m
\end{equation}
for all $m\in M$, where, if $({}_im)_{i=\ov{1, n}}$ is a basis of $M$ 
with dual basis
$({}^im)_{i=\ov{1, n}}$ in $M^*$, then $ {}^{*i}m$ is the image of 
${}_im$ under the
canonical map $M\to M^{**}$.
Moreover, the morphism $\theta'_M$ is defined
by the same formula (\ref{thetay}) as $\theta_M$. A straightforward but 
tedious computation proves that the morphism $\Theta _M:\ M^*\ra {}^*M$
defined by (\ref{rl}) is given by
\begin{eqnarray}
&&\hspace*{-2cm}\Theta _M(m^*)=
\langle m^*, S(p^1)f^2\cd (g^1\cd {}_jm)_{(0)}\rangle \nonumber\\
&&\hspace*{5mm}
\langle {}^jm, S(q^2)\smi (q^1\smi (f^1(g^1\cd {}_jm)_{(- 
1)}g^2)p^2)\cd {}_im\rangle {}^im,
\label{rly}
\end{eqnarray}
for all $m^*\in M^*$. The inverse map $\Theta _M^{-1}:\ {}^*M\ra M^*$ 
is given by (see (\ref{lr}))
\begin{equation}
\Theta ^{-1}_M({}^*m)=
\langle {}^*m, V^2\cd {}_jm_{(0)}\rangle \langle {}^jm, \b 
S(V^1{}_jm_{(-1)})\cd {}_im\rangle {}^im,\label{irly}
\end{equation}
for any ${}^*m\in {}^*M$, where
\begin{equation}\label{v}
V=V^1\ot V^2=\smi (f^2p^2)\ot \smi (f^1p^1).
\end{equation}
Therefore, we obtain the following result.

\begin{proposition}\label{prop1}
Let $H$ be a quasi-Hopf algebra and $M, N\in \yd ^{\rm fd}$. Then
$^r\Gamma _M={\theta'}^{-1}_M\circ \Theta _{M^*}:\ M^{**}\ra M$ is
an isomorphism of Yetter-Drinfeld modules. Explicitly, ${}^r\Gamma _M$ 
is given by
\begin{equation}\label{gr}
^r\Gamma _M(m^{**})=\langle m^{**}, {}^im\rangle
q^1S^{-2}((S^2(p^1)\cd {}_im)_{(-1)}^{\mathcal F})p^2S^2(q^2)\cd
(S^2(p^1)\cd {}_im)_{(0)}^{\mathcal F}
\end{equation}
for all $m^{**}\in M^{**}$, where
\begin{equation}\label{calf}
{\mathcal F}={\mathcal F}^1\ot {\mathcal F}^2
=S(g^2)f^1\ot S(g^1)f^2,
\end{equation}
and where, in general, if ${\mf F}={\mf F}^1\ot {\mf F}^2$ 
is a twist on $H$ and $M\in \yd$ then we denote 
$$
\l _M^{\mf F}(m)=m_{(- 1)}^{\mf F}\ot m_{(0)}^{\mf F}:=
{\mf F}^1({\mf G}^1\cd m)_{(-1)}{\mf G}^2\ot {\mf F}^2\cd
({\mf G}^1\cd m)_{(0)},
$$
where ${\mf G}={\mf G}^1\ot {\mf G}^2={\mf F}^{-1}$. 

The inverse of $^r\Gamma _M$ is given by
$^r\Gamma _M^{-1}=\Theta ^{-1}_{M^*}\circ \theta '_M$, that is,
\begin{equation}\label{igr}
^r\Gamma ^{-1}_M(m)=\langle {}^im, S((p^1\cd m)_{(-1)}p^2)\a \cd
(p^1\cd m)_{(0)}\rangle {}^{*i}m.
\end{equation}
Similarly, ${}^l\Gamma _M:=\theta ^{-1}_M\circ \Theta ^{-1}_{{}^*M}:\ 
{}^{**}M\ra M$
provides an isomorphism of Yetter-Drinfeld modules. Explicitly, we have 
that
\begin{equation}\label{gl}
^l\Gamma_{M}({}^{**}m)=\langle {}^{**}m, 
{}^im\rangle 
V^1g^1S((\smi (V^2g^2)\cd {}_im)_{(-1)})\a \cd (\smi (V^2g^2)\cd 
{}_im)_{(0)},
\end{equation}
for all ${}^{**}m\in {}^{**}M$. The inverse of ${}^l\Gamma _M$ is given 
by
\begin{equation}\label{igl}
^l\Gamma ^{-1}_M(m)=
\langle {}^im, S^{-2}(q^1(p^1\cd m)_{(-1)}p^2)q^2\cd (p^1\cd 
m)_{(0)}\rangle {}^{*i}m.
\end{equation}
\end{proposition}

\begin{proof}
Straightforward, we leave the details to the reader.
\end{proof}

Now, we will focus on the second isomorphism.
Let ${\mathcal C}$ be a rigid braided monoidal category and $M, N\in {\mathcal 
C}$.
By \cite{ag} there exist natural isomorphisms
$$\sigma ^*_{M, N}:\ M^*\ot N^*\ra (M\ot N)^*$$
$${}^*\sigma _{M, N}:\ {}^*M\ot {}^*N\ra {}^*(M\ot N).$$
In fact, $\sigma^*_{M, N}=\phi ^*_{N, M}\circ c^{-1}_{N^*, M^*}$, where
$\phi_{N, M}^*:\ N^*\ot M^*\ra (M\ot N)^*$ is the composition
\begin{eqnarray}
&&\hspace*{-15mm}
 \phi_{N, M}^*=\ev_N\ot {(M\ot N)^*}
~~\circ~~  ({N^*}\ot (\ev_M\ot N))\ot {(M\ot N)^*}\nonumber\\
&\circ& {N^*}\ot a^{-1}_{M^*, M, N}\ot {(M\ot N)^*}
~~\circ~~ a_{N^*, M^*, M\ot N}\ot {(M\ot N)^*}\nonumber\\
&\circ& a^{-1}_{N^*\ot M^*, M\ot N, (M\ot N)^*}
~~\circ~~{N^*\ot M^*}\ot \coev_{M\ot N}\label{phir}
\end{eqnarray}
with inverse
\begin{eqnarray}
&&\hspace*{-13mm}
\phi ^{*-1}_{N, M}= (\ev_{M\ot N}\ot {N^*})\ot {M^*}
~~\circ~~ a^{-1}_{(M\ot N)^*, M\ot N, N^*}\ot {M^*}\nonumber\\
&&\circ~~ ({(M\ot N)^*}\ot a^{-1}_{M, N, N^*})\ot {M^*}
~~\circ~~ ({(M\ot N)^*}\ot (M\ot \coev_N))\ot {M^*}\nonumber\\
&&\circ~~ 
 a^{-1}_{(M\ot N)^*, M, M^*}
~~\circ ~~{(M\ot N)^*}\ot \coev_M.
\label{iphir}
\end{eqnarray}
In a similar way, we can define ${}^*\phi _{N, M}:\ {}^*N\ot {}^*M\ra 
{}^*(M\ot N)$ and its 
inverse ${}^*\phi_{N, M}^{-1}$; then we put 
${}^*\sigma _{M, N}={}^*\phi _{N, M}\circ c^{-1}_{{}^*N, 
{}^*M}$.

\begin{proposition}\label{prop2}
Let $H$ be a quasi-Hopf algebra and $M, N\in \yd^{\rm fd}$.
Consider $\{{}_im\}_{i=\ov{1, s}}$ and $\{{}^im\}_{i=\ov{1, s}}$ be dual 
bases in $M$ and $M^*$, and $\{{}_jn\}_{j=\ov{1, t}}$ and 
${\{}^jn\}_{j=\ov{1, t}}$ 
dual bases in $N$ and $N^*$. Then the map 
$\sigma ^*_{M, N}:\ M^*\ot N^*\ra (M\ot N)^*$ defined by
\begin{eqnarray}
&&\sigma ^*_{M,N}(m^*\ot n^*)(m\ot n)\nonumber\\
&&\hspace*{1cm}
=\langle m^*, f^2\tilde{q}^2_2X^3\smi (\tilde{q}^1X^1(p^1\cd 
n)_{(-1)}p^2)\cd m\rangle 
\langle n^*, f^1\tilde{q}^2_1X^2\cd (p^1\cd 
n)_{(0)}\rangle \label{ydsr}
\end{eqnarray}
is an isomorphism
of Yetter-Drinfeld modules.  Here $p_R=p^1\ot p^2$ and
$q_L=\tilde{q}^1\ot \tilde{q}^2$ are the elements defined in 
(\ref{qr}) and (\ref{qla}). The inverse of $\sigma ^*_{M, N}$ is given by
\begin{equation}\label{iydsr}
\sigma ^{*-1}_{M, N}(\mu )=
\langle \mu , (g^1\cd {}_jn)_{(-1)}g^2\cd {}_im\ot (g^1\cd 
{}_jn)_{(0)}\rangle {}^im\ot {}^jn.
\end{equation}
In a similar way, the map
${}^*\sigma _{M, N}:\ {}^*M\ot {}^*N\ra {}^*(M\ot N)$ defined by
\begin{eqnarray}
&&{}^*\sigma _{M, N}({}^*m\ot {}^*n)(m\ot n)
=\langle {}^*m, \smi (f^2q^1(x^1\tilde{p}^1\smi (f^1)\cd 
n)_{(-1)}x^2\tilde{p}^2_1)\cd m\rangle \nonumber\\
&&\hspace*{3cm}
\langle {}^*n, \smi (x^3\tilde{p}^2_2)q^2\cd (x^1\tilde{p}^1\smi 
(f^1)\cd n)_{(0)}\rangle
\label{ydsl}\end{eqnarray}
is an isomorphism of Yetter-Drinfeld modules; $q_R=q^1\ot q^2$ and $p_L=
\tilde{p}^1\ot \tilde{p}^2$ are the elements defined in 
(\ref{qra}) and (\ref{ql}). 
The inverse of ${}^*\sigma _{M, N}$ is given by
\begin{equation}\label{iydsl}
{}^*\sigma _{M, N}^{-1}(\nu )=
\langle \nu , (\smi (g^2)\cd {}_jn)_{(-1)}\smi (g^1)\cd {}_im 
\ot (\smi (g^2)\cd {}_jn)_{(0)}\rangle {}^im\ot {}^jn.
\end{equation}
\end{proposition}

\begin{proof}
We will show that $\sigma^*_{M, N}=\phi^*_{N, M}\circ c_{N^*, 
M^*}^{-1}$, where
$\phi^*_{N, M}$ is given by (\ref{phir}).
We first calculate $\phi ^*_{N,M}$.
As before, we write $q_R=q^1\ot q^2$ and $p_R=p^1\ot p^2=P^1\ot P^2$, 
and then compute for all $n^*\in N^*$, $m^*\in M^*$, $m\in M$ and
$n\in N$ that:
\begin{eqnarray*}
&&\hspace*{-2cm}
\phi ^*_{N, M}(n^*\ot m^*)(m\ot n)\\
&{{\rm (\ref{phir})}\atop =}&\langle m^*, S(x^1X^2y^1_2)\a x^2(X^3y^2\b 
S(y^3))_1\cd m\rangle \\
&&\hspace*{1cm}
\langle n^*, S(X^1y^1_1)\a x^3(X^3y^2\b S(y^3))_2\cd 
n\rangle \\
&{{\rm (\ref{gd})}\atop =}&\langle m^*, S(y^1_2)\gamma ^1(y^2\b 
S(y^3))_1\cd m\rangle
\langle n^*, S(y^1_1)\gamma ^2(y^2\b S(y^3))_2\cd n\rangle \\
&{{\rm (\ref{ca}, \ref{qr})}\atop =}&\langle m^*, f^1S(p^1)_1g^1\gamma 
^1p^2_1\cd m\rangle
\langle n^*, f^2S(p^1)_2g^2\gamma ^2p^2_2\cd n\rangle \\
&{{\rm (\ref{gdf})}\atop =}&\langle m^*, f^1(S(p^1)\a p^2)_1\cd m\rangle
\langle n^*, f^2(S(p^1)\a p^2)_2\cd n\rangle \\
&{{\rm (\ref{qr}, \ref{q6})}\atop =}&\langle m^*, f^1\cd m\rangle \langle 
n^*, f^2\cd n\rangle .
\end{eqnarray*}
Using (\ref{iphir}) or by a direct computation it is easy to see that
\begin{equation}\label{sat}
\phi ^{*-1}_{N, M}(\mu )=
\langle \mu , g^1\cd {}_im\ot g^2\cd {}_jn\rangle
{}^jn\ot {}^im
\end{equation}
for all $\mu \in (M\ot N)^*$. Also, it is not hard to see
that (\ref{g2}) and (\ref{pf}) imply that
\begin{equation}\label{ufo}
q^1g^1_1\ot S(q^2g^1_2)g^2=S(X^3)f^1\ot S(X^1\b S(X^2)f^2).
\end{equation}
The same relations and (\ref{l3a}) imply that
\begin{equation}\label{ufox}
S(p^2f^1)F^1f^2_1\ot S(p^1)F^2f^2_2=q_L.
\end{equation}
Using (\ref{ufo}), the axioms of a quasi-Hopf algebra,
again (\ref{g2}) and (\ref{pf}), and
finally (\ref{pql}) we obtain the following relation:
\begin{eqnarray}
&&q^1[g^1S(\tilde{q}^2_1X^2)f^2]_1\ot
S(q^2[g^1S(\tilde{q}^2_1X^2)f^2]_2)g^2S(\tilde{q}^1X^1)
\ot S(\tilde{q}^2X^3)f^1\nonumber\\
&&\hspace*{3cm}=f^2F^2_2x^2\ot S(F^2x^3)\ot f^1F^1_1x^1.\label{uf}
\end{eqnarray}
Let $f=F^1\ot F^2={\bf F}^1\ot {\bf F}^2$
and $p_R=P^1\ot P^2$ be other copies of $f$ and $p_R$.
For $m^*\in M^*$, $n^*\in N^*$, $m\in M$ and $n\in N$
we calculate:
\begin{eqnarray*}
&&\hspace*{-2cm}
(\phi ^*_{N, M}\circ c^{-1}_{N^*, M^*})(m^*\ot n^*)(m\ot n)\\
&{{\rm (\ref{y6}, \ref{rdy2})}\atop =}&\langle m^*, \tilde{q}^1X^1
\smi (F^1(g^1S(\tilde{q}^2_1X^2)f^2\cd n)_{(-1)}g^2)
p^2S(\tilde{q}^2_2X^3)f^1\cd m\rangle \\
&&\langle n^*, S(p^1)F^2\cd (g^1S(\tilde{q}^2_1X^2)f^2\cd 
n)_{(0)}\rangle \\
&{{\rm (\ref{y3p}, \ref{uf})}\atop =}&\langle m^*, \smi ({\bf F}^1f^2_1
F^1_{(2, 1)}x^2_1(P^1\cd n)_{(-1)}P^2S(F^2x^3))
p^2f^1F^1_1x^1\cd m\rangle \\
&&\langle n^*, S(p^1){\bf F}^2f^2_2F^1_{(2, 2)}x^2_2\cd (P^1\cd 
n)_{(0)}\rangle \\
&{{\rm (\ref{uf}, \ref{ufox})}\atop =}&\langle m^*, \smi 
(\tilde{q}^1F^1_{(2, 1)}x^2_1
(P^1\cd n)_{(-1)}P^2S(F^2x^3))F^1_1x^1\rangle \\
&&\langle n^*, F^1_2\tilde{q}^2x^2_2\cd (P^1\cd n)_{(0)}\rangle \\
&{{\rm (\ref{ql1a}, \ref{q3})}\atop =}&\langle m^*, F^2\tilde{q}^2_2X^3
\smi (\tilde{q}^1X^1(P^1\cd n)_{(-1)}P^2)\cd m\rangle \\
&&
\langle n^*, F^1\tilde{q}^2_1X^2\cd (P^1\cd n)_{(0)}\rangle \\
&{{\rm (\ref{ydsr})}\atop =}&\sigma ^*_{M, N}(m^*\ot n^*)(m\ot n).
\end{eqnarray*}
Obviously, the inverse of $\sigma^*_{M, N}$ is
$\sigma ^{*-1}_{M, N}=c_{N^*, M^*}\circ \phi ^{*-1}_{N, M}$. By 
(\ref{sat}) it 
follows that $\sigma^{*-1}_{M, N}$ is defined by the formula given in 
(\ref{iydsr}).
The assertion concerning the morphism ${}^*\sigma _{M, N}$ can be 
proved in a
similar way. We only notice that
$$
{}^*\phi _{N, M}({}^*n\ot {}^*m)(m\ot n)=
\langle {}^*m, \smi (f^2)\cd m\rangle \langle {}^*n, \smi (f^1)\cd 
n\rangle
$$
for all ${}^*m\in {}^*M$, ${}^*n\in {}^*N$, $m\in M$ and $n\in N$.
\end{proof}

Let $(H, R)$ be a quasitriangular quasi-Hopf algebra. Then we have a 
monoidal fuctor
$\mathbb{F}:\  {}_H{\mathcal M}\ra \yd $ which acts as the identity on 
objects and morphisms. For 
$M\in {}_H{\mathcal M}$, $\mathbb{F}(M)=M$ as a left $H$-module, and 
with left $H$-coaction given by
\begin{equation}\label{funct}
\l _M(m)=R^2\ot R^1\cd m,
\end{equation}
for all $m\in M$. Moreover, this functor sends algebras, coalgebras, 
bialgebras etc.
in ${}_H{\mathcal M}$ to the corresponding objects in
$\yd $ (see \cite[Proposition 2.4]{bn}).

\begin{corollary}(\cite[Lemma 3.2]{bpv2}).\label{co1}
Let $H$ be a quasitriangular quasi-Hopf algebra and
$M$ a finite dimensional left $H$-module. Then
$M\cong M^{**}$ and $M\cong {}^{**}M$ as $H$-modules.
\end{corollary}

\begin{proof}
We have seen above that $M$ can be viewed as an object in $\yd $, so 
$M\cong M^{**}$ and
$M\cong {}^{**}M$ as Yetter-Drinfeld modules, cf. Proposition 
\ref{prop1}. Thus,
$M\cong M^{**}$ and $M\cong {}^{**}M$ as $H$-modules. We will write 
down explicitly these isomorphisms. We have that ${}^r\Gamma _M: M^{**}\ra M$ 
is given by 
$$
{}^r\Gamma _{M}(m^{**})=
\langle m^{**}, {}^im\rangle u^{-1}\cd {}_im,
$$
with inverse 
$$
{}^r\Gamma _M^{-1}(m)=\langle {}^im, u\cd m\rangle 
{}^{*i}m,
$$
for all $m^{**}\in M^{**}$ and $m\in M$, which is equivalent to
${}^r\Gamma _M^{-1}(m)(\varphi )=\varphi (u\cd m)$ for all
$\varphi \in M^*$ and $m\in M$. In this way we recover the 
isomorphism
$\hat{\psi }:\ M\ra M^{**}$ given in \cite{bpv2}.

Similarly, ${}^l\Gamma _M:\ {}^{**}M\ra M$ is given by
$$
^l\Gamma ({}^{**}m)=\langle {}^{**}m, 
{}^im\rangle u\cd {}_im
$$
for all ${}^{**}m\in {}^{**}M$, with inverse 
$$
{}^l\Gamma _M^{-1}(m)=\langle {}^im, u^{-1}\cd m\rangle {}^{*i}m
$$ 
for all $m\in M$.
\end{proof}

\begin{remark}\label{re1}
The element $u$ plays a central role in the theory of
quasitriangular Hopf algebras. It is therefore natural to try to 
generalize its
properties to the quasi-Hopf algebra setting. The major problem
in \cite{bn3} was to find the inverse of the element $u$ defined in 
\cite{ac}. Now, if we forget
the definitions of $u$ and $u^{-1}$ in the quasi-Hopf case, combining 
Corollary \ref{co1} with the
similar result in the Hopf case we will obtain in a natural way the 
definitions for $u$ and its inverse $u^{-1}$.
\end{remark}

\begin{corollary}(\cite[Lemma 3.3]{bpv2}).\label{co2}
Let $(H, R)$ be a quasitriangular quasi-Hopf algebra and $M$, $N$ two 
finite dimensional
left $H$-modules. Then $M^*\ot N^*\cong (M\ot N)^*$ and ${}^*M\ot 
{}^*N\cong {}^*(M\ot N)$
as $H$-modules.
\end{corollary}

\begin{proof}
We regard $M, N$ as objects in $\yd ^{\rm fd}$ via the functor 
$\mathbb{F}$ defined in (\ref{funct}).
Thus, by Proposition \ref{prop2}, we obtain that $M^*\ot N^*\cong (M\ot 
N)^*$ and
${}^*M\ot {}^*N\cong {}^*(M\ot N)$ as Yetter-Drinfeld modules, so they 
are isomorphic
also as $H$-modules. Moreover, (\ref{invr2}) implies that the 
isomorphism
$\sigma ^*_{M, N}:\ M^*\ot N^*\ra (M\ot N)^*$ defined in (\ref{ydsr}) 
is given by
$$
\sigma ^*_{M, N}(m^*\ot n^*)(m\ot n)=\langle m^*, f^2\ov {R}^2\cd 
m\rangle \langle n^*,
f^1\ov {R}^1\cd n\rangle
$$
for all $m^*\in M^*$, $n^*\in N^*$, $m\in M$, $n\in N$, where
$R^{-1}:=\sum \ov{R}^1\ot \ov {R}^2$ is the inverse of the $R$-matrix 
$R$. Note that
it is just the isomorphism $\mu _{M, N}$ defined in \cite{bpv2}. Also, 
(\ref{invr1})
implies that the isomorphism ${}^*\sigma _{M, N}: {}^*M\ot {}^*N\ra 
{}^*(M\ot N)$
defined by (\ref{ydsl}) is given by
$$
{}^*\sigma _{M, N}({}^*m\ot {}^*n)(m\ot n)=
\langle {}^*m, \smi (f^2\ov {R}^2)\cd m\rangle \langle {}^*n, \smi 
(f^1\ov{R}^1)\cd n\rangle ,
$$
for all ${}^*m\in {}^*M$, ${}^*n\in {}^*N$, $m\in M$ and $n\in N$. 
Finally, it is not hard to see that
$$
\sigma ^{*-1}_{M, N}(\mu )=
\langle \mu , R^2g^2\cd {}_im\ot R^1g^1\cd {}_jn\rangle {}^im\ot {}^jn
$$
for all $\mu \in (M\ot N)^*$, and
\begin{eqnarray*}
{}^*\sigma ^{-1}_{M, N}(\nu )
&{{\rm (\ref{ydsl}, \ref{funct})}\atop =}&
\langle \nu , R^2\smi (g^1)\cd {}_im\ot R^1\smi (g^2)\cd {}_jn\rangle 
{}^im\ot {}^jn\\
&{{\rm (\ref{ext})}\atop =}&
\langle \nu , \smi (R^2g^2)\cd {}_im\ot \smi (R^1g^1)\cd {}_jn\rangle 
{}^im\ot {}^jn
\end{eqnarray*}
for any $\nu \in {}^*(M\ot N)$, completing the proof.
\end{proof}

\begin{remark}
Continuing the ideas of Remark \ref{re1}, we notice that the above 
Corollary
suggests the two formulae (\ref{invr1}) and (\ref{invr2}) for the 
inverse
of the $R$-matrix $R$, and, also, the formula $(\ref{ext})$. All these
formulae where first proved by Hauser and Nill \cite{hn2} in the case
that $(H, R)$ is a finite dimensional quasitriangular quasi-Hopf 
algebra. They
also used the bijectivity of the antipode.
\end{remark}

\section{Applications}\selabel{5}
\setcounter{equation}{0}
Let $B$ be a Hopf algebra in a braided rigid 
category $\Cc$. Then 
$B^*$ and ${}^*B$ are also Hopf algebras in $\Cc$, see \cite{ag,ta}.
The structure maps on $B^*$ are the following ones:
\begin{eqnarray}
&&m_{B^*}:\ B^*\ot B^*\rTo^{\sigma ^*_{B, B}}(B\ot B)^*
\rTo^{\Delta _B^*}B^*,\label{mbra}\\
&&\Delta _{B^*}:\ B^*\rTo^{m_B^*}(B\ot B)^*
\rTo^{\sigma ^{*-1}_{B, B}}B^*\ot B^*,\label{combra}\\
&&S_{B^*}=S_B^*,\mbox{${\;}$}
u_{B^*}=\va _B^*, \mbox{${\;}$}
\va _{B^*}=u_B^*.\label{rebra}
\end{eqnarray}
The structure maps on the left dual can be obtained in a similar way.\\
We can consider algebras, coalgebras, bialgebras and Hopf algebras in 
the
braided category $\yd$. More precisely, an algebra in $\yd$ is an
object $B\in \yd$ with the additional structure:
\begin{itemize}
\item[-] $B$ is a left $H$-module algebra, i.e.
$B$ has a multiplication
and a unit $1_B$ satisfying the
conditions:
\begin{eqnarray}
(ab)c&=&(X^1\cd a)[(X^2\cd b)(X^3\cd c)],\nonumber\\
h\cd (ab)&=&(h_1\cd a)(h_2\cd b),\nonumber\\
h\cd 1_B&=&\va (h)1_B,\label{mal}
\end{eqnarray}
for all $a, b, c\in B$ and $h\in H$.
\item[-] $B$ is a quasi-comodule
algebra, that is, the multiplication $m$ and the unit
$\eta $ of $B$ intertwine the $H$-coaction $\lambda _B$. By (\ref{y4}),
this means:
\begin{eqnarray}
&&\hspace*{-2cm}
\lambda _B(bb')=X^1(x^1Y^1\cd b)_{(-1)}x^2(Y^2\cd b')_{(-1)}
Y^3 \nonumber\\
&\ot& [X^2\cd (x^1Y^1\cd b)_{(0)}]
[X^3x^3\cd (Y^2\cd b')_{(0)}],
\label{qca1}
\end{eqnarray}
for all $b, b'\in B$, and
\begin{equation}\label{qca2}
\lambda _B(1_B)=1_H\ot 1_B.
\end{equation}
\end{itemize}
A coalgebra in $\yd$ is an object $B$ with
\begin{itemize}
\item[-] $B$ is a left $H$-module coalgebra, i.e. $B$ has a 
comultiplication
$\un {\Delta }_B :\ B\ra B\ot B$ (we will denote
$\un {\Delta }(b)=b_{\un {1}}\ot b_{\un {2}}$)
and a usual counit $\un {\va }_B$ such that:
\begin{eqnarray}
&&X^1\cd b_{(\un{1}, \un{1})}
\ot X^2\cd b_{(\un{1}, \un{2})}
\ot X^3\cd b_{\un{2}}=b_{\un{1}}\ot b_{(\un{2}, \un{1})}\ot
b_{(\un{2}, \un{2})},\label{mc1}\\
&&
\un {\Delta }_B(h\cd b)=h_1\cd b_{\un{1}}\ot h_2\cd b_{\un{2}},
\mbox{${\;\;\;}$}
\un{\va }_B(h\cd b)=\va (h)\un{\va }_B(b),\label{mc2}
\end{eqnarray}
for all $h\in H$, $b\in B$, where we adopt
for the quasi-coassociativity
of $\un {\Delta }_B$  the same notations as in
the \seref{1}.
\item[-] $B$ is a quasi-comodule coalgebra,
i.e. the comultiplication
$\un {\Delta }_B$ and the counit $\un {\va }_B$ intertwine the
$H$-coaction $\l _B$. Explicitly, for all $b\in B$ we must have that:
\begin{eqnarray}
&&b_{(-1)}\ot b_{(0)_{\un{1}}}\ot b_{(0)_{\un{2}}}=
X^1(x^1Y^1\cd b_{\un{1}})_{(-1)}
x^2(Y^2\cd b_{\un{2}})_{(-1)}Y^3\ot\nonumber\\
&&\hspace*{1cm}X^2\cd (x^1Y^1\cd b_{\un{1}})_{(0)}\ot
X^3x^3\cd (Y^2\cd b_{\un{2}})_{(0)},\label{qcc1}
\end{eqnarray}
and
\begin{equation}
\un{\va }_B(b_{(0)})b_{(-1)}=\un{\va }_B(b)1.\label{qcc2}
\end{equation}
\end{itemize}
An object $B\in \yd$ is a bialgebra if it is an algebra and a coalgebra 
in $\yd $ such that
$\un {\Delta }_B$ is an algebra morphism, i.e.
$\un {\Delta }_B(1_B)=1_B\ot 1_B$ and, by (\ref{bi})
and (\ref{y5}), for all $b, b'\in B$ we have that:
\begin{eqnarray}
&&\hspace*{-2cm}\Delta _B(bb')=[y^1X^1\cd b_{\un{1}}]
[y^2Y^1(x^1X^2\cd b_{\un{2}})_{(-1)}x^2
X^3_1\cd b'_{\un{1}}]\nonumber\\
&&\ot [y^3_1Y^2\cd (x^1X^2\cd b_{\un{2}})_{(0)}]
[y^3_2Y^3x^3X^3_2\cd b'_{\un{2}}].\label{by}
\end{eqnarray}
Finally, a bialgebra $B$ in $\yd $ is a Hopf algebra
if there exists a morphism $\un {S}: B\ra B$ in $\yd $ such that
$\un{S}(b_{\un{1}})b_{\un{2}}=b_{\un{1}}S(b_{\un{2}})=
\un{\va }_B(b)1$ for all $b\in B$; we then say that
$B$ is a braided Hopf algebra.\\
If $B\in \yd ^{\rm fd}$ is a Hopf algebra then $B^*$ and ${}^*B$
become Hopf algebras in $\yd $. We will compute the structure maps.

\begin{proposition}\prlabel{5.1}
Let $H$ be a quasi-Hopf algebra and $B$ a finite dimensional Hopf 
algebra in $\yd $.
Then $B^*$ is a Hopf algebra in $\yd $ with the following structure 
maps:
\begin{eqnarray}
\hspace*{-1.6cm}
(a^*\un{\circ}\hspace*{1mm}b^*)(b)&=&\langle a^*, 
f^2\tilde{q}^2_2X^3\smi (\tilde{q}^1X^1
(p^1\cd b_{\un{2}})_{(-1)}p^2)\cd b_{\un{1}}\rangle \nonumber\\
&&\langle b^*, f^1\tilde{q}^2_1X^2\cd (p^1\cd
b_{\un{2}})_{(0)}\rangle \\
\Delta _{B^*}(b^*)&=&
\langle b^*, [(g^1\cd {}_jb)_{(-1)}g^2\cd {}_ib](g^1\cd 
{}_jb)_{(0)}\rangle
{}^ib\ot {}^jb\\
S_{B^*}(b^*)&=&b^*\circ \un{S},
\mbox{${\;\;}$}
u_{B^*}(1)=\un{\va }_B,
\mbox{${\;\;}$}
\va _{B^*}(b^*)=b^*(1_B),
\end{eqnarray}
for all $a^*, b^*\in B^*$ and $b\in B$, where $\{{}_ib\}_{i=\ov{1, n}}$ 
is a basis in $B$
with dual basis $\{{}^ib\}_{i=\ov{1, n}}$ in $B^*$. Similarly, ${}^*B$ 
is a Hopf
algebra in $\yd $ with the following structure maps:
\begin{eqnarray}
\hspace*{1cm}
({}^*a\hspace*{1mm}\ov{\circ}\hspace*{1mm}{}^*b)(b)&=&
\langle {}^*a, \smi (f^2q^1(x^1\tilde{p}^1\smi (f^1)\cd 
b_{\un{2}})_{(-1)}
x^2\tilde{p}^2_1)\cd b_{\un{1}}\rangle \nonumber\\
&&\langle {}^*b, \smi (x^3\tilde{p}^2_2)q^2\cd (x^1\tilde{p}^1\smi 
(f^1)\cd b_{\un{2}})_{(0)}\rangle \\
\Delta _{{}^*B}(b^*)&=&
\langle {}^*b, [(\smi(g^2)\cd {}_jb)_{(-1)}
\smi (g^1)\cd {}_ib](\smi (g^2)\cd {}_jb)_{(0)}\rangle 
{}^ib\ot {}^jb\\
S_{{}^*B}({}^*b)&=&{}^*b\circ \un{S},
\mbox{${\;\;}$}
u_{{}^*B}(1)=\un{\va }_B,
\mbox{${\;\;}$}
\va _{{}^*B}({}^*b)={}^*b(1_B),
\end{eqnarray}
for all ${}^*a, {}^*b\in {}^*B$ and $b\in B$.
\end{proposition}

\begin{proof}
Follows easily by computing the morphisms in (\ref{mbra}-\ref{rebra}), 
using
(\ref{ydsr}-\ref{iydsl}); the details are left to the reader.
\end{proof}

Let $H$ be a Hopf algebra. It is well-known that $H$ becomes an algebra
in the monoidal category $\yd $, with Yetter-Drinfeld structure given by
$$
h\triangleright h'=h_1h'S(h_2),
\mbox{${\;\;\;}$}
\lambda (h)=h_1\otimes h_2
$$
for all $h, h'\in H$. Moreover, $H$ is quantum commutative as an algebra
in $\yd $ (see for example \cite{cvoz}). Now, let $H$ be a
quasi-Hopf algebra. In \cite{bpv} a new multiplication $\circ$
on $H$ is introduced, given by the formula
\begin{equation}\label{ma}
h\circ h'=X^1hS(x^1X^2)\a
x^2X^3_1h'S(x^3X^3_2)
\end{equation}
Let $H_0$ be equal to $H$ as a vector space, with multiplication 
$\circ$.
Then $H_0$
is a left $H$-module algebra with unit $\beta $ and left $H$-action
given by
\begin{equation}\label{s1}
h\tr h'=h_1h'S(h_2),
\end{equation}
for all $h, h'\in H$. It was also shown in \cite{bn} that $H_0$
is an algebra in the category $\yd $, with $H$-coaction given by
\begin{equation}\label{s2}
\lambda _H(h)=h_{(-1)}\ot h_{(0)}=
X^1Y^1_1h_1g^1S(q^2Y^2_2)Y^3\ot X^2Y^1_2h_2g^2S(X^3q^1Y^2_1),
\end{equation}
where $f^{-1}=g^1\otimes g^2$ and $q=q_R=q^1
\otimes q^2$ are the elements defined by (\ref{g}) and
(\ref{qra}). In addition, in \cite{bcp}
it was shown that $H_0$ is actually quantum commutative as an algebra in
$\yd $.\\
If $(H, R)$ is quasitriangular, then $H_0$ is a Hopf algebra with 
bijective antipode 
in ${}_H{\mathcal M}$, with the additional structure (see \cite{bn}):
\begin{eqnarray}
&&
\un{\Delta }(h)=h_{\un{1}}\ot h_{\un{2}}\nonumber\\
&&\hspace*{8mm}=x^1X^1h_1g^1S(x^2R^2y^3X^3_2)
\ot x^3R^1\tr 
y^1X^2h_2g^2S(y^2X^3_1),\label{und}\\
&&\un{\va }(h)=\va (h),\label{unva}\\
&&\un{S}(h)=X^1R^2p^2S(q^1(X^2R^1p^1\tr h)S(q^2)X^3),\label{unant}
\end{eqnarray}
for all $h\in H$, where $R=R^1\ot R^2$ and
$f^{-1}=g^1\ot g^2$, $p_R=p^1\ot p^2$ and $q_R=q^1\ot 
q^2$
are the elements defined by (\ref{g}), (\ref{qr}) and (\ref{qra}).
If we consider the left $H$-coaction
\begin{equation}
\lambda _{H_0}(h)=R^2\ot R^1\tr h,\label{scshz}\\
\end{equation}
induced by (\ref{funct}), then $H_0$ becomes a Hopf algebra in $\yd $,
with bijective antipode. 
From now on, we will refer to $H_0$ as a Hopf algebra in 
${}_H{\mathcal M}$, and, via the monoidal functor 
${\mathbb{F}}:\ {}_H{\mathcal M}\to \yd $, in $\yd$, 
with structure maps (\ref{ma}), (\ref{s1}) and 
(\ref{und}-\ref{scshz}).

If $H$ is finite dimensional, then $H_0^*$ is also
a Hopf algebra in $\yd $. By \thref{3.2} and (\ref{ext}), 
$H_0^*$ is a Yetter-Drinfeld module via
\begin{eqnarray}
&&\hspace*{-2cm}
(h\mapsto \v )(h')=\v (S(h)\tr h')\\
&&\hspace*{-2cm}
\l _{H_0^*}(\v )=R^2\ot R^1\mapsto \v
\end{eqnarray}
for all $\v \in H^*$ and $h, h'\in H$. The structure of $H_0^*$ as a
Hopf algebra in $\yd $ is given by:
\begin{eqnarray}
(\v \un{\circ} \Psi )(h)&=&\langle \v , f^2\ov{R}^2\tr 
h_{\un{1}}\rangle
\langle \Psi , f^1\ov{R}^1\tr h_{\un{2}}\rangle ,\label{dmhz1}\\
1_{H_0^*}&=&\un{\va },\label{duhz}\\
\Delta _{H_0^*}(\v )&=&
\langle \v , (R^2g^2\tr {}_ie)\circ (R^1g^1\tr {}_je)\rangle
{}^ie\ot {}^je\label{dcmhz},\\
\va _{H_0^*}(\v )&=&\v (\b ),\label{dchz}\\
S_{H_0^*}(\v )&=&\v \circ \un{S},\label{danthz}
\end{eqnarray}
for all $h\in H$ and $\v , \Psi \in H^*$, where $R^{-1}=\sum 
\ov{R}^1\ot \ov{R}^2$,
$\{{}_ie\}_{i=\ov{1, n}}$ is a basis of $H$ and
$\{{}^ie\}_{i=\ov{1, n}}$ the corresponding dual basis of $H^*$. The 
left dual
${}^*H_0$ of $H_0$ is also a Hopf algebra
in $\yd $. First, by \thref{3.2} and (\ref{ext}),
${}^*H_0$ is a Yetter-Drinfeld module via
\begin{eqnarray}
&&\hspace*{-2cm}
(h\succ \v )(h')=\v (\smi (h)\tr h')\\
&&\hspace*{-2cm}
\l _{{}^*H_0}(\v )=R^2\ot R^1\succ \v
\end{eqnarray}
for all $\v \in H^*$ and $h, h'\in H$. Then the structure of ${}^*H_0$ 
as a Hopf algebra in $\yd $ is given by the formulae
\begin{eqnarray}
&&(\v \ov {\circ}\Psi )(h)=\langle \v , \smi (f^2\ov{R}^2)\tr 
h_{\un{1}}\rangle
\langle \Psi , \smi (f^1\ov{R}^1)\tr h_{\un{2}}\rangle ,\label{dmhz2}\\
&&1_{{}^*H_0}=\un{\va},
\end{eqnarray}
\begin{eqnarray}
&&\Delta _{{}^*H_0}(\v )=
\langle \v , [\smi (R^2g^2)\tr {}_ie]\circ [\smi (R^1g^1)\tr 
{}_je]\rangle {}^ie\ot {}^je,\label{dcmhz2}\\
&&\va _{{}^*H_0}(\v )=\v (\b ),\label{dchz2}\\
&&S_{{}^*H_0}(\v )=\v \circ \un{S},
\end{eqnarray}
for all $\v \in H^*$ and $h, h'\in H$. From (\ref{rl}) we know
that $\Theta _{H_0}:\ H_0^*\ra {}^*H_0$ is an isomorphism
of Yetter-Drinfeld modules. In this particular case, $\Theta_{H_0}$ is 
given by
\begin{eqnarray*}
\Theta _{H_0}(\v )&{{\rm (\ref{rly}, \ref{funct})}\atop 
=}&\langle \v , S(p^1)f^2R^1g^1\tr {}_je\rangle \\
&&\hspace*{1cm}\langle {}^je, S(q^2)\smi (q^1\smi (f^1R^2g^2)p^2)\tr 
{}_ie\rangle {}^ie\\
&{{\rm (\ref{ext})}\atop =}&
\langle \v , S(R^1p^1)\tr {}_je\rangle
\langle {}^je, \smi (q^1R^2p^2S^2(q^2))\tr {}_ie\rangle {}^ie\\
&{{\rm (\ref{inelmu})}\atop =}&\langle \v , \smi 
(u^{-1})\tr {}_ie\rangle {}^ie
=u^{-1}\succ \v
\end{eqnarray*}
for all $\v \in H_0^*$. Since $S^2(u)=u$ it is easy to see that $\Theta 
_{H_0}^{-1}(\v )=
u\mapsto \v $ for all $\v \in {}^*H_0$.

\begin{proposition}\prlabel{5.2}
Let $(H, R)$ be a triangular quasi-Hopf algebra. With notation as above,
$\Theta_{H_0}:\ H_0^*\ra {}^*H_0$ is a braided Hopf algebra
isomorphism.
\end{proposition}
\begin{proof}
It is well known that in a symmetric monoidal category ${\cal C}$ the 
canonical isomorphism $\Theta _B$, where $B$ is a Hopf algebra in $\Cc$, 
provides a braided Hopf algebra isomorphism between $B^*$ and ${}^*B$. 
Since $(H, R)$ is triangular, the category ${}_H{\cal M}^{\rm fd}$ is 
symmetric, and this finishes the proof.    
\end{proof}

\begin{remark}\relabel{5.3}
If $(H, R)$ is not triangular then $\Theta _{H_0}$ is, in general, not
an algebra or a coalgebra morphism. Indeed, if $(H, R)$ is an arbitrary
quasitriangular quasi-Hopf algebra then 
computations similar to the ones presented above show that
\begin{eqnarray*}
&&\Theta _{H_0}((\ov{r}^1\ov{R}^2\mapsto \v )
\un{\circ}(\ov{r}^2\ov{R}^1\mapsto \Psi ))=\Theta _{H_0}(\v )\ov{\circ}
\Theta _{H_0}(\Psi ),\\
&&(\Delta _{{}^*H_0}\Theta _{H_0})(\v )=R_{21}R\succ
(\Theta _{H_0}\ot \Theta _{H_0})(\Delta _{H_0^*}(\v ))
\end{eqnarray*}
for any $\v , \Psi \in H_0^*$, where $\ov{r}^1\ot \ov{r}^2$ is 
another copy of $R^{-1}$ 
and we extend the action of $H$ on ${}^*H_0$ to 
an action of $H\ot H$ on ${}^*H_0\ot {}^*H_0$.
\end{remark}

Let $H$ be a quasi-Hopf algebra. Then $H^*$ is an $(H,H)$-bimodule, by
$$
\langle h\rh \v , h'\rangle =\v (h'h), \mbox{${\;\;\;}$}
\langle \v \lh h, h'\rangle =\v (hh').
$$
The convolution $\langle \v \Psi , h\rangle =\v (h_1)\psi (h_2)$, 
$h\in H$, 
is a multiplication on $H^*$; it is not associative, 
but only quasi-associative: 
\begin{equation}\label{mbia1}
(\v \psi)\xi=(X^1\rh \v \lh x^1)[(X^2\rh \psi \lh x^2)
(X^3\rh \xi \lh x^3)], \mbox{${\;\;\;}$$\forall \v , \psi , \xi
\in H^*$.}
\end{equation}
In addition, for all $h\in H$ and $\v , \psi \in H^*$ we have
that
\begin{equation}\label{mbia2}
h\rh (\v \psi )=(h_1\rh \v )(h_2\rh \psi )
\mbox{${\;\;\;}$and${\;\;\;}$}
(\v \psi )\lh h=(\v \lh h_1)(\psi \lh h_2).
\end{equation}
By \cite{bc1}, if $H$ is a finite dimensional quasitriangular
quasi-Hopf algebra then on
the dual of $H$ there exists another structure of Hopf algebra in $\yd 
$,
denoted by $\un {H}^*$. The structure of $\un{H}^*$ as a
Yetter-Drinfeld module is given
by the formulae
\begin{eqnarray}
h\cd \v &=&h_1\rh \v \lh \smi (h_2),\label{mhs}\\
\l _{\un{H}^*}(\v )&=&R^2\ot R^1\cd \v ,\label{chs}
\end{eqnarray}
for all $h\in H$, $\varphi \in H^*$.
The structure of $\un {H}^*$ as a Hopf algebra in $\yd$ is given by:
\begin{eqnarray}
&&\v \bullet \Psi =(x^1X^1\rh \v \lh \smi 
(f^2x^3_2Y^3R^1X^2))\nonumber\\
&&\hspace*{1cm}
(x^2Y^1R^2_1X^3_1\rh \Psi \lh \smi 
(f^1x^3_1Y^2R^2_2X^3_2)),\label{mbdd}\\
&&\Delta _{\un{H}^*}(\v )=X^1_1p^1\rh \v _2\lh \smi (X^1_2p^2)\ot 
X^2\rh \v _1\lh \smi (X^3),\label{cbdd}\\
&&\va _{\un{H}^*}(\v )=\v (\smi (\a )),\label{cbdd1}\\
&&S_{\un{H}^*}(\v )=Q^1q^1R^2x^2\cdot [p^1P^2S(Q^2) \nonumber\\
&&\hspace*{2cm}
\rh \ov {S}^{-1}(\v )\lh S(q^2R^1x^1P^1)x^3\smi (p^2)],\label{anbdd}
\end{eqnarray}
for all $\v , \Psi \in H^*$, where  $Q^1\ot Q^2$ is another copy
of $q_R$ and $\ov{S}^{-1}(\v )=\v \circ \smi $ for any $\v \in H^*$.
The unit element is $\va$.

\begin{proposition}\prlabel{5.4}
Let $(H, R)$ be a finite dimensional triangular quasi-Hopf algebra.
Then ${}^*H_0$ and $\un{H}^{*cop}$ are isomorphic as braided Hopf 
algebras.
\end{proposition}

\begin{proof}
Since $(H, R)$ is triangular it follows that $c^{-1}_{N, M}=c_{M, N}$, 
so $\un{H}^{*cop}$ is a Hopf algebra in $\yd$. On the other hand, 
both ${}^*H_0$ and $\un{H}^{*cop}$ are Hopf algebras in $\yd $ as 
images, through
the functor $\mathbb{F}$ defined by (\ref{funct}), of corresponding 
objects in ${}_H{\mathcal M}$.
Therefore, it suffices to
prove that ${}^*H_0$ and $\un{H}^{*cop}$ are isomorphic as Hopf 
algebras in
${}_H{\mathcal M}$. To this end we claim that $\mu :\ {}^*H_0\ra 
\un{H}^{*cop}$ given by
$$
\mu (\v )=g^1\rh \v \lh \smi (g^2)
$$
is a Hopf algebra isomorphism in ${}_H{\mathcal M}$. In fact, $\mu $ is 
$H$-linear since
\begin{eqnarray*}
\mu (h\succ \v )(h')&=&\langle \v , \smi (h)\tr (\smi 
(g^2)h'g^1)\rangle \\
&{{\rm (\ref{ca})}\atop =}&\langle \v , \smi (h_2g^2)h'h_1g^1\rangle \\
&{{\rm (\ref{mhs})}\atop =}&\langle \mu (\v ), \smi (h_2)h'h_1\rangle =(h\cd 
\mu (\v ))(h')
\end{eqnarray*}
for all $\v \in H^*$ and $h, h'\in H$. It is also an algebra morphism in
${}_H{\mathcal M}$ since
\begin{eqnarray*}
&&\hspace*{-3cm}\mu (\v \ov{\circ}\Psi )(h)\\
&{{\rm (\ref{dmhz2}, \ref{und}, \ref{ca})}\atop =}&\langle \v , \smi 
(f^2\ov{R}^2)\tr x^1X^1\smi (F^2g^2_2{\bf G}^2)h_1g^1_1G^1\\
&&\hspace*{-2cm} S(x^2R^2y^3X^3_2)\rangle
\langle \Psi , \smi (f^1\ov{R}^1)x^3R^1\tr y^1X^2\\
&&\hspace*{-2cm}
\smi (F^1g^2_1{\bf G}^1)h_2g^1_2G^2S(y^2X^3_1)\rangle \\
&{{\rm (\ref{ca}, \ref{g2}, \ref{pf})}\atop =}&\langle \v , \smi 
(f^2\ov{R}^2)\tr x^1\smi (F^2X^3g^2)
h_1X^1_1g^1_{(1, 1)}{\bf G}^1_1G^1\\
&&\hspace*{-2cm}S(x^2R^2y^3)\rangle
\langle \Psi , \smi (f^1\ov{R}^1)x^3R^1\tr y^1\\
&&\hspace*{-2cm}\smi (F^1X^2g^1_2{\bf G}^2)h_2
X^1_2g^1_{(1, 2)}{\bf G}^1_2G^2S(y^2)\rangle \\
&{{\rm (\ref{g2}, \ref{pf}, \ref{q1})}\atop =}&
\langle \v , \smi (f^2\ov{R}^2)\tr x^1\smi (F^2X^3g^2)h_1X^1_1y^1g^1_1\\
&&\hspace*{-2cm}{\bf G}^1S(x^2R^2)\rangle
\langle \Psi , \smi (f^1\ov{R}^1)\tr \smi (F^1X^2y^3g^1_{(2, 2)}\\
&&\hspace*{-2cm}{\bf G}^2_2G^2S(x^3_1R^1_1))
h_2X^1_2y^2g^1_{(2, 1)}
{\bf G}^2_1G^1S(x^3_2R^1_2)\rangle \\
&{{\rm (\ref{ca}, \ref{ext}, \ref{qt3})}\atop =}&
\langle \v , \smi (f^2\ov{R}^2)\tr x^1\smi (F^2X^3g^2)
h_1X^1_1y^1R^2g^1_2{\bf G}^2S(x^2)\rangle \\
&&\hspace*{-2cm}
\langle \Psi , \smi (f^1\ov{R}^1)\tr 
\smi (F^1X^2y^3R^1_2(g^1_1{\bf G}^1S(x^3))_2G^2)\\
&&\hspace*{-2cm}h_2
X^1_2y^2R^1_1(g^1_1{\bf G}^1S(x^3))_1G^1\rangle 
\end{eqnarray*}
\begin{eqnarray*}
&{{\rm (\ref{g2}, \ref{pf})}\atop =}&\langle \v , 
\smi (f^2\ov{R}^2)\tr 
\smi (F^2X^3x^3g^2_2{\bf G}^2)h_1X^1_1y^1R^2\\
&&\hspace*{-2cm}x^2g^2_1{\bf G}^1\rangle
\langle \Psi , \smi (f^1\ov{R}^1)\tr \smi 
(F^1X^2y^3R^1_2x^1_2g^1_2G^2)h_2X^1_2y^2R^1_1x^1_1g^1_1G^1\rangle \\
&{{\rm (\ref{ca}\;\;twice)}\atop =}&
\langle \v , \smi (F^2X^3x^3\ov{R}^2_2g^2)h_1X^1_1y^1R^2x^2
\ov{R}^2_1g^1\rangle \\
&&\hspace*{-2cm}\langle \Psi , \smi (R^1x^1\ov{R}^1)\tr \smi 
(F^1X^2y^3G^2)h_2X^1_2y^2G^1\rangle \\
&{{\rm (\ref{qt2}, \ref{q3})}\atop =}&\langle \v , \smi 
(F^2x^3_2X^3\ov{R}^2Y^2g^2)h_1x^1Y^1g^1\rangle \\
&&\hspace*{-2cm}\langle \Psi , \smi (\ov{R}^1Y^3)\tr \smi 
(F^1x^3_1X^2G^2)h_2x^2X^1G^1\rangle \\
&{{\rm (\ref{ca},\;\;R^{-1}=R_{21})}\atop =}&\langle \mu (\v ), \smi 
(F^2x^3_2X^3R^1Y^2)h_1x^1Y^1\rangle \\
&&\hspace*{-2cm}\langle \mu (\Psi ), \smi 
(F^1x^3_1X^2R^2_2Y^3_2)h_2x^2X^1R^2_1Y^3_1\rangle \\
&{{\rm (\ref{mbdd})}\atop =}&(\mu (\v )\bullet \mu (\Psi ))(h)
\end{eqnarray*}
for all $\v , \Psi \in H^*$, $h\in H$, and
$\mu (1_{{}^*H_0})=\mu (\va )=\va =1_{\un{H}^{*cop}}$. It remains
to show that $\mu $ is a coalgebra morphism and a bijection. To this 
end, we calculate for any $\v \in H^*$:
\begin{eqnarray*}
&&\hspace*{-3cm}
(\mu \ot \mu )\Delta _{{}^*H_0}
=\langle \v , (\smi (R^2g^2)\tr {}_ie)\circ
(\smi (R^1g^1)\tr {}_je)\rangle \\
&&G^1\rh {}^ie\lh \smi (G^2)
\ot {\bf G}^1\rh {}^je\lh \smi ({\bf G}^2)\\
&{{\rm (\ref{ma}, \ref{ext})}\atop =}&
\langle \v , X^1[R^2\smi (g^1)\tr \smi (G^2){}_ieG^1]S(x^1X^2)\a 
x^2X^3_1\\
&&[R^1\smi (g^2)\tr
\smi ({\bf G}^2){}_je{\bf G}^1]S(x^3X^3_2)\rangle {}^ie\ot {}^je\\
&{{\rm (\ref{ca})}\atop =}&
\langle \v , X^1R^2_1\smi (g^1_2G^2){}_ieg^1_1G^1S(x^1X^2R^2_2)\a x^2\\
&&[X^3R^1\smi (g^2)\tr
\smi ({\bf G}^2){}_je{\bf G}^1]S(x^3)\rangle {}^ie\ot {}^je\\
&{{\rm (\ref{qt2}, \ref{ext}, \ref{qt3})}\atop =}&
\langle \v , X^1\smi (y^2r^2g^2_2G^2){}_iey^1g^1S(x^1R^2X^3)\a x^2\\
&&[R^1X^2\smi (y^3r^1g^2_1G^1)
\tr \smi ({\bf G}^2){}_je{\bf G}^1]S(x^3)\rangle {}^ie\ot {}^je\\
&{{\rm (\ref{g2}, \ref{pf})}\atop =}&
\langle \v , \smi (y^2r^2X^3g^2){}_iey^1X^1g^1_1G^1S(x^1R^2)\a x^2\\
&&[R^1\smi (y^3r^1X^2g^1_2G^2)\tr \smi ({\bf G}^2){}_je{\bf 
G}^1]S(x^3)\rangle {}^ie\ot {}^je\\
&{{\rm (\ref{ext}, \ref{qt3}, \ref{ca})}\atop =}&
\langle \v , \smi (R^2_2X^3g^2){}_ieR^2_1X^2g^1_2G^2S(x^1)\a x^2\\
&&\smi (R^1_2X^1_2g^1_{(1, 2)}G^1_2{\bf G}^2)
{}_jeR^1_1X^1_1g^1_{(1, 1)}G^1_1{\bf G}^1S(x^3)
\rangle {}^ie\ot {}^je\\
&{{\rm (\ref{g2}, \ref{pf}, \ref{q1})}\atop =}&
\langle \v , \smi (R^2_2X^3g^2){}_ieR^2_1X^2x^3
g^1_{(2, 2)}G^2_2{\bf G}^2\a \\
&&\smi (R^1_2X^1_2x^2g^1_{(2, 1)}G^2_1{\bf 
G}^1){}_jeR^1_1X^1_1x^1g^1_1G^1S(x^3)\rangle
{}^ie\ot {}^je\\
&{{\rm (\ref{l3a}, \ref{q5})}\atop =}&
\langle \v , \smi (R^2_2X^3g^2){}_ieR^2_1X^2x^3\smi (R^1_2X^1_2x^2\b )\\
&&{}_jeR^1_1X^1_1x^1g^1\rangle {}^ie\ot {}^je\\
&{{(\ref{qr}, \ref{mhs})}\atop =}&
\langle \mu (\v ), \smi (X^3){}_ieX^2\smi 
(X^1_2p^2){}_jeX^1_1p^1\rangle \\
&&R^2\cd {}^ie\ot R^1\cd {}^je\\
&{{\rm (\ref{cbdd},\;\;R^{-1}=R_{21})}\atop =}&
c^{-1}_{\un{H}^*, \un{H}^*}(\Delta _{\un{H}^*}(\mu (\v )))
=\Delta _{\un{H}^{*cop}}(\mu (\v ))
\end{eqnarray*}
as needed. Obviously $\va \mu =\va $. It is easy to see that $\mu $ is 
bijective with inverse
$$
\mu ^{-1}(\v )=f^1\rh \v \lh \smi (f^2)
$$
for any $\v \in H^*$. Thus, the proof is complete.
\end{proof}

\begin{corollary}\colabel{5.5}
If $(H, R)$ is a finite dimensional triangular quasi-Hopf algebra then
$H^*_0\cong {}^*H_0\cong \un{H}^{*cop}$ as braided Hopf algebras.
\end{corollary}

\begin{center}
{\sc Acknowledgement}
\end{center}
The authors wish to thank the referee for his suggestion to use the centre
construction to simplify the proofs of the results in Sections \ref{se:2}
and \ref{se:3}.


\end{document}